\documentclass[reqno]{amsart}
\usepackage{amscd,amsmath,amssymb,amsfonts,amsthm,verbatim}
\usepackage{times}
\usepackage[cmtip, all]{xy}
\usepackage{graphicx}
\usepackage[active]{srcltx}
\usepackage{color}
\usepackage{tikz-cd}
\usetikzlibrary{calc}
\usepackage{hyperref}

\usepackage[color, notref, notcite]{}     % refs and labels
\definecolor{refkey}{gray}{.5}   % graylevel for refs
\definecolor{labelkey}{gray}{.5} % graylevel for labels
\definecolor{Red}{rgb}{1,0,0}

\newcommand{\ra}{\rightarrow}

\newcommand{\pf}{{\bf Proof : }}

\newcommand{\qedwhite}{\hfill \ensuremath{\Box}}

\newtheorem{theo}{Theorem}[section]
\newtheorem{prop}[theo]{Proposition}
\newtheorem{lem}[theo]{Lemma}
\newtheorem{defi}[theo]{Definition}
\newtheorem{cor}[theo]{Corollary}

\newtheorem{notn}[theo]{Notation}
\newtheorem{rem}[theo]{Remark}

\newcommand{\ga}{\alpha}

\newcommand{\Spec}{\mbox{\rm Spec\,}}

\title{A note on relative Vaserstein symbol}
\author{ Kuntal Chakraborty}

\newcommand{\Addresses}{{% additional braces for segregating \footnotesize
  \bigskip
  \footnotesize

 \textsc{Kuntal Chakraborty, School of Mathematics, Tata Institute of Fundamental Research,\\ \noindent
           1, Dr. Homi Bhabha Road, Mumbai 400005, INDIA
   } \par\nopagebreak
  \textit{E-mail}: Kuntal Chakraborty \texttt{<kuntal@math.tifr.res.in>}

  \medskip

  }}

\begin{document}
\maketitle

\section{\bf{Introduction}}
In an unpublished work of Fasel-Rao-Swan in \cite{fr} the notion of the relative Witt group $W_E(R,I)$ is defined. In this article we will give the details of this construction:

Let $I$ be an ideal of $R$. We denote the set \{$ \alpha \in Alt^{\prime}_{2n}(R): \alpha\equiv \chi_r^{\pm 1}\perp \chi_{n-r}^{\pm 1} \pmod I~ \text{for some r} \geq 0$\}  by $Alt^{\prime}_{2n}(R,I)$, where $Alt^{\prime}_{2n}(R)$ is the set of all skew-symmetric matrices in $GL_{2n}(R)$ and $\chi_1=\begin{pmatrix}
0&1\\
-1&0
\end{pmatrix}$ and $\chi_{r+1}=\chi_r\perp \chi_1$ for $r\geq 1$. For $m\leq n$, $Alt^{\prime}_{2m}(R,I)$ can be embedded into $Alt^{\prime}_{2n}(R,I)$ by the inclusion map $i_{m,n}:Alt^{\prime}_{2m}(R,I)\hookrightarrow Alt^{\prime}_{2n}(R,I)$ given by $ \alpha\longmapsto \alpha\perp\chi_{n-m} $. We put $ Alt^{\prime}(R,I)=\varinjlim Alt^{\prime}_{2n}(R,I) $.
We can define an equivalence relation $\sim_I$ on $Alt^{\prime}(R,I)$ as follows:

 Let $ \alpha \in Alt^{\prime }_{2m}(R,I) $ and $\beta \in Alt^{\prime}_{2n}(R,I)$. Then $\alpha \sim_I \beta$ if and only if there exists $t\in \mathbb{N}$ and $\epsilon \in E_{2(m+n+t)}(R,I)$ such that :
 \begin{center}
 $\alpha \perp \chi_{n+t}=\epsilon^T(\beta \perp \chi_{m+t})\epsilon$
 \end{center} 
  We denote the set $Alt^{\prime}(R,I)/\sim_I$ by $W^{\prime}_E(R,I)$.

  We can also define the set $Alt_{2n}(R,I)$ as the subset of $Alt^{\prime}_{2n}(R,I)$ consisting of skew-symmetric matrices of Pfaffian $1$. We put $Alt(R,I):= \varinjlim Alt_{2n}(R,I)$. As  the relation $\sim_I$ defines an equivalence relation on $Alt^{\prime}(R,I)$, then $\sim_I$ defines an equivalence relation on $Alt(R,I)$. We denote the set $Alt(R,I)/\sim_I$ by $W_E(R,I)$. 

  This definition is compatible with the definition of the symplectic Witt group introduced by Vaserstein  in \cite{SV}. This definition can also be realised as a kernel of a certain map between two Witt groups, namely the following sequence is exact: 

\[
\begin{tikzcd}
0\arrow{r}{} &W_E(R,I)\arrow{r}{i} &W_E(D)\arrow{r}{(p_1)_*} &W_E(R)\arrow{r}{}&0
\end{tikzcd}
\]
where $D$ is the double ring of $R$ with respect to the ideal $I$ and $p_1$ is the first projection map.

Likewise in \cite{SV}, we can consider the relative Vaserstein symbol , namely, the map $$V_{R,I}: Um_3(R,I)/E_3(R,I)\rightarrow W_E(R,I).$$

In \cite{SV}, Vaserstein proved that for a commutative ring $R$ of dimension $2$, the Vaserstein symbol $V_R: Um_3(R)/E_3(R)\rightarrow W_E(R)$ is bijective.

 In \cite{rv}, Rao-van der Kallen proved that the Vaserstein symbol $V_R: Um_3(R)/E_3(R)\rightarrow W_E(R)$ is bijective for $R$ smooth affine algebra of dimension $3$ over a $C_1$- field $k$ which is perfect if its characteristic is $2$ or $3$.

Also in \cite{frs}, Fasel-Rao-Swan proved that the Vaserstein symbol $V_R:Um_3(R)/E_3(R)\rightarrow W_E(R)$ is injective for $R$ smooth affine algebra of dimension $4$ over an algebraically closed field.

It is believed that these results can be generalised in relative cases. In fact we have:

\begin{theo}
$($See Corollary \ref{5.10}$)$
Let $R$ be a commutative ring of dimension $2$ and $I\subset R$ be an ideal of $R$. Then the Vaserstein symbol $V_{R,I}: Um_3(R,I)/E_3(R,I)\rightarrow W_E(R,I)$ is bijective. 
\end{theo}
 
 \begin{theo}
 $($See Theorem \ref{localcompinte}$)$
 Let $R$ be an affine non-singular algebra of dimension $3$ over a perfect $C_1$-field $k$ and $I\subset R$ be a local complete intersection ideal of $R$. Then the Vaserstein symbol $V_{R,I}: Um_3(R,I)/E_3(R,I)\rightarrow W_E(R,I)$ is bijective.
 \end{theo}

  \begin{theo}
  $($See Theorem \ref{localcomintedim4}$)$
 Let $R$ be an affine non-singular algebra of dimension $4$ over an algebraically closed field $k$ and $I\subset R$ be a local complete intersection ideal of $R$. Then the Vaserstein symbol $V_{R,I}: Um_3(R,I)/E_3(R,I)\rightarrow W_E(R,I)$ is injective.
 \end{theo}
 
 In fact we studied more about injectivity of Vaserstein symbol. We have considered two cases:  One injectivity of the Vaserstein symbol $V_{R[X], I[X]}$ where $R$ is an affine algebra of dimension $\leq 3$ over a field $k$ and $I$ is an ideal of $R$. The next is injectivity of the Vaserstein symbol $V_{R[X],(X^2-X)}$ where $R$ is an affine algebra of dimension $\leq 3$ over a perfect $C_1$-field $k$. In fact we have:
 
 \begin{theo}
$($See Theorem \ref{extendedvas}$)$
Let $R$ be a non-singular affine algebra of dimension $2$ over a perfect $C_1$-field $k$, and let $I\subset R$ be a local complete intersection ideal of $R$. Then the Vaserstein symbol $$V_{R[X],I[X]}:MSE_3(R[X],I[X])\rightarrow W_E(R[X],I[X]))$$ is injective.
\end{theo}

\begin{theo}
$($see Theorem \ref{relvas}$)$
Let $R$ be an affine algebra of dimension $2$ over a perfect $C_1$-field $k$ with $6R=R$. Then the Vaserstein symbol $V_{R[X],(X^2-X)}: Um_3(R[X],(X^2-X))/E_3(R[X],(X^2-X))\rightarrow W_E(R[X],(X^2-X))$ is bijective.
\end{theo}
 We study the above case because of Fasel-Rao-Swan. According to Fasel-Rao-Swan the relative Witt group $W_E(R[X],(X^2-X))$ of a two dimensional non-singular algebra $R$ over a perfect $C_1$-field was a divisible group prime to characteristic. And divisibility of $W_E(R[X],(X^2-X))$ gives completability of a relative unimodular row of size $3$ over $R[X]$ via a result of Suslin.

Later we computed the kernel of the Vaserstein map $V_{R,I}$: 

\begin{theo}
$($see Theorem \ref{kernel}$)$
Let $R$ be a ring of dimension $3$ and $I\subset R$ be an ideal of $R$. If the orbit space $MSE_4(R,I)$ has a nice group structure, then,  the Vaserstein map $V_{R,I}$ induces a bijection  of the map $$\phi: Um_3(R,I)/\{\sigma\in SL_3(R,I)\cap E_5(R,I)\}\equiv W_E(R,I).$$
\end{theo}

It is believed that for a $3$-dimensional affine algebra non-singularity is not necessary for establishing injectivity of the Vaserstein symbol. We will give an example to show this at the end of this article. In \cite{rs}, it is showed that  a Rees algebra $R[It]$, with $R$ non-singular is non-singular if and only if $I=0,R$ and $I$ is generated by a single element. By this result we will be able to construct a singular $3$-dimensional algebra over a perfect $C_1$-field for which the Vaserstein symbol is injective. Moreover we can extend this result for $4$-dimensional singular algebra over a field $k$. In fact we have:

\begin{theo}
$($ See Theorem \ref{4.10}$)$
Let $R$ be a non-singular affine algebra of dimension $3$ over a field $k$ and $I\subset R$ be an ideal. Then the Vaserstein symbol $V_{R[It]}: MSE_3(R[It])\rightarrow W_E(R[It])$ is injective.
\end{theo}
 We established the injectivity of Vaserstein symbol over extended Rees algebra also:
 
 \begin{theo}
 $($ See Theorem \ref{extdedRees}$)$
Let $R$ be a ring of dimension $2$ and $I\subset R$ be an ideal of $R$. Then the Vaserstein symbol $V_{R[It,t^{-1}]}: MSE_3(R[It,t^{-1}])\rightarrow W_E(R[It,t^{-1}])$ is bijective.
\end{theo}

\section{\bf{Witt Group and Relative Witt Group}}
           
 $R$ will denote a commutative ring with $1\neq 0$, unless stated otherwise. $I$ will denote an ideal of $R$.

We denote by $GL_n(R,I)$ by the kernel of the canonical map $GL_{n}(R) \longrightarrow GL_{n}\left(\frac{R}{I}\right)$. Let $SL_{n}(R,I)$ denotes the subgroup of $GL_{n}(R,I)$ of elements of determinant $1$.

\begin{defi}
${\bf{The ~Relative ~Elementary ~group ~E_n(R,I):} }$
Let $R$ be a ring and $I\subset R$ be an ideal. The relative elementary group is the subgroup of $SL_{n}(R,I)$ generated by the matrices of the form $\alpha e_{i,j}(a) \alpha^{-1}$, where $\alpha \in E_{n}(R),i\neq j$ and $a\in I$.
\end{defi}

We identify $GL_{n}(R)$ with a subgroup of $GL_{n+1}(R)$ by associating the matrix $\begin{pmatrix}
1 &0\\
0 &\alpha
\end{pmatrix}\in GL_{n+1}(R)$, with the element $\alpha \in GL_{n}(R)$.
We set \begin{center}
$GL(R)=\varinjlim GL_{n}(R),  GL(R,I)=\varinjlim GL_{n}(R,I)$
\\$SL(R)=\varinjlim SL_{n}(R),  SL(R,I)=\varinjlim SL_{n}(R,I)$
\\$E(R)=\varinjlim E_{n}(R),  E(R,I)=\varinjlim E_{n}(R,I)$
\end{center}
Let $Alt^{\prime}_{2n}(R)$ be the set of all skew-symmetric matrices in $GL_{2n}(R)$. For any $r\in \mathbb{N}$, define $\chi_{r}$ inductively by 
\begin{center}
$\chi_{1}=\begin{pmatrix}
0&1\\
-1&0
\end{pmatrix}$
\end{center}
and $\chi_{r+1}=\chi_{r}\perp\chi_{1}$. Clearly $\chi_{r} \in Alt^{\prime}_{2r}(R)$. For any $m\leq n$ , $Alt^{\prime}_{2m}(R)$ can be embedded into $Alt^{\prime}_{2n}(R)$ by the inclusion map $i_{m,n}:Alt^{\prime}_{2m}(R) \hookrightarrow Alt^{\prime}_{2n}(R)$ given by $\alpha \mapsto \alpha\perp\chi_{n-m}$. Define $Alt^{\prime}(R)=\varinjlim Alt^{\prime}_{2n}(R)$. For $\alpha \in Alt^{\prime}_{2m}(R)$ and $\beta\in Alt^{\prime}_{2n}(R)$ define $\alpha\sim_E\beta$ if and only if there exists some $t \in \mathbb{N}$ and $\varepsilon \in E_{2(m+n+t)}(R)$ such that $\alpha \perp \chi_{n+t}=\varepsilon^{T}(\beta \perp \chi_{m+t})\varepsilon$. It can be shown that this relation is reflexive, symmetric and transitive i.e. $\sim_E$ is an equivalence relation on $Alt^{\prime}(R)$. It can be shown that $Alt^{\prime}(R)/\sim_E$ is an abelian group with respect to the operation $\perp$\cite[Section 3]{SV}. Denote this group by $W^{\prime}_E(R)$. Similarly one can consider the subset $Alt_{2n}(R)\subset Alt^{\prime}_{2n}(R)$, the set of all skew-symmetric matrices in $GL_{2n}(R)$ with Pfaffian $1$. Using the same embedding, as defined earlier, one can construct the set $Alt(R):=\varinjlim Alt_{2n}(R)$. Also $\sim_E$ is an equivalence relation on $Alt(R)$ and $Alt(R)/\sim_E$ is an abelian group with respect to the operation $\perp$. Denote this group by $W_E(R)$. This group is called \textit{ the elementary symplectic Witt group}.

Let $I$ be an ideal of $R$. We denote the set \{$ \alpha \in Alt^{\prime}_{2n}(R): \alpha\equiv \chi_r^{\pm 1}\perp \chi_{n-r}^{\pm 1} \pmod I~ \text{for some r} \geq 0$\}  by $Alt^{\prime}_{2n}(R,I)$. For $m\leq n$, $Alt^{\prime}_{2m}(R,I)$ can be embedded into $Alt^{\prime}_{2n}(R,I)$ by the inclusion map $i_{m,n}:Alt^{\prime}_{2m}(R,I)\hookrightarrow Alt^{\prime}_{2n}(R,I)$ given by $ \alpha\longmapsto \alpha\perp\chi_{n-m} $. We put $ Alt^{\prime}(R,I)=\varinjlim Alt^{\prime}_{2n}(R,I) $.
We can define an equivalence relation $\sim_I$ on $Alt^{\prime}(R,I)$ as follows:

 Let $ \alpha \in Alt^{\prime }_{2m}(R,I) $ and $\beta \in Alt^{\prime}_{2n}(R,I)$. Then $\alpha \sim_I \beta$ if and only if there exists $t\in \mathbb{N}$ and $\epsilon \in E_{2(m+n+t)}(R,I)$ such that :
 \begin{center}
 $\alpha \perp \chi_{n+t}=\epsilon^T(\beta \perp \chi_{m+t})\epsilon$
 \end{center} 
  We denote the set $Alt^{\prime}(R,I)/\sim_I$ by $W^{\prime}_E(R,I)$.

  We can also define the set $Alt_{2n}(R,I)$ as the subset of $Alt^{\prime}_{2n}(R,I)$ consisting of skew-symmetric matrices of Pfaffian $1$. We put $Alt(R,I):= \varinjlim Alt_{2n}(R,I)$. As  the relation $\sim_I$ defines an equivalence relation on $Alt^{\prime}(R,I)$, then $\sim_I$ defines an equivalence relation on $Alt(R,I)$. We denote the set $Alt(R,I)/\sim_I$ by $W_E(R,I)$. 
  
  \begin{defi}
 Let $f_1:A\rightarrow C$  be a ring homomorphism between two commutative rings $A$ and $C$. Also let $f_2:B\rightarrow C$ be another ring homomorphism between $B$ and $C$. A ring $D$ is said to be fibre product of $A$ and $B$ with respect to $C$ if the following two conditions are satisfied:
 
 1) there exist two ring homomorphisms,namely $p_1:D\rightarrow A$, $p_2:D\rightarrow B$ such that the following diagram is commutative:
 \[
\begin{tikzcd}
D \arrow{r}{p_1} \arrow{d}{p_2}
               &A\arrow{d}{f_1}\\
B \arrow{r}{f_2} &C
\end{tikzcd}
\]

2) If there exists some ring E with two ring homomorphisms $q_1:E\rightarrow A$, $q_2:E\rightarrow B$ such that the following diagram is commutative :
\[
\begin{tikzcd}
E \arrow{r}{q_1} \arrow{d}{q_2}
               &A \arrow{d}{f_1}\\
B \arrow{r}{f_2} &C
\end{tikzcd}
\]
 then there exists a unique map $h:E\rightarrow D$ such that the following diagram is commutative:
\[
\begin{tikzcd}
E\arrow[bend left]{drr}{q_1}
\arrow[bend right]{ddr}[swap]{q_2}
\arrow[dotted]{dr}[description]{\exists ! h} & & \\
& D \arrow{r}{p_1} \arrow{d}{p_2} & A \arrow{d}{f_1} \\
& B \arrow{r}{f_2} & C
\end{tikzcd}
\]
We denote the fibre product of $A$ and $B$ with respect to $C$ by $A \times_C B$.
\end{defi}

 Consider the fibre product
\[
\begin{tikzcd}
D \arrow{r}{p_1} \arrow{d}{p_2}
               &R \arrow{d}{\pi}\\
R \arrow{r}{\pi} &R/I
\end{tikzcd}
\]
Here the ring $D$ is called the double ring of the ring $R$ with respect to the ideal $I$. $D$ can be identified with the set \{ $(a,b)\in R\times R:a-b\in I$\}

\begin{lem}
\label{1.5}
Let $\alpha\in M_n(D)$. Then $\alpha$ can be associated canonically to an element $(\alpha_1,\alpha_2)\in M_n(R)\times M_n(R)$ such that  $\alpha_1 \equiv \alpha_2 \pmod{I} $. In other words $M_n(D)$ is isomorphic to the double ring of $M_n(R)$ with itself with respect to the ideal $M_n(I)$.
\end{lem}

${\pf}$ The maps $p_1$ and  $p_2$ induces the ring homomorphisms $p_1^*:M_n(D)\rightarrow M_n(R)$ and  $p_2^* : M_n(D)\rightarrow M_n(R)$. We have a commutative diagram 
\[
\begin{tikzcd}
M_n(D)\arrow{r}{(p_1)^*}\arrow{d}{(p_2)^*} &M_n(R)\arrow{d}{\pi^*}\\
M_n(R)\arrow{r}{\pi^*} &M_n(R/I)
\end{tikzcd}
\]

Since the fibre product is unique upto isomorphism  it is enough to show that $M_n(D)\cong M_n(R)\times_{M_n(R/I)}M_n(R)$.
By the universal property of  the fibre product there exists a unique map $h:M_n(D)\rightarrow M_n(R)\times_{M_n(R/I)}M_n(R)$,which is defined by $h(A)=(A_1,A_2)$ where $A=(a_{ij})$ with $a_{ij}=(a^1_{ij},a^2_{ij})$ such that $a^1_{ij}-a^2_{ij}\in I$, $A_1=(a^1_{ij})$ and $A_2=(a^2_{ij})$.
Now define $g:M_n(R)\times_{M_n(R/I)}M_n(R)\rightarrow M_n(D)$ as following:

Let $(B_1,B_2)\in M_n(R)\times_{M_n(R/I)}M_n(R)$. Then $\bar{B_1}=\bar{B_2}$ in $M_n(R/I)$. If $B_1=(b^1_{ij})$ and $B_2=(b^2_{ij})$ then we have $b^1_{ij}-b^2_{ij}\in I$. Therefore the matrix $B \in M_n(D)$ is well-defined. We define $g(B_1,B_2)=B$. It can be shown that $g$ and $h$ both are ring homomorphisms. Now we have $g\circ h=1_{M_n(D)}$ and $h\circ g=1_{M_n(R)\times_{M_n(R/I)}M_n(R)}$. Hence $M_n(D)\cong M_n(R)\times_{M_n(R/I)}M_n(R)$.

 $~~~~~~~~~~~~~~~~~~~~~~~~~~~~~~~~~~~~~~~~~~~~~~~~~~~~~~
\qedwhite$

Using the above canonical isomorphism we always represent the element of $M_n(D)$ as an element of $M_n(R)\times_{M_n(R/I)} M_n(R)$, i.e. as a pair of two matrices satisfying certain properties.

\begin{lem}
Let $M,N\in M_n(D)$ be such that $M=(M_1,M_2)$ and $N=(N_1,N_2)$. Then $M\perp N=(M_1\perp N_1,M_2\perp N_2)$ in $M_{2n}(D)$.
\end{lem}
 
 $\pf$ This is true in $M_{2n}(R)\times M_{2n}(R)$ hence the relation is true in $M_{2n}(R)\times_{M_{2n}(R/I)} M_{2n}(R)$. Finally the relation is true in $M_{2n}(D)$ by Lemma \ref{1.5}.
 
 $~~~~~~~~~~~~~~~~~~~~~~~~~~~~~~~~~~~~~~~~~~~~~~~~~~~
 \qedwhite$
 
\begin{defi}
Let $R$ be a ring and $I\subset R$ be an ideal. The excision ring of $R$ with respect to the ideal $I$ is denoted by $R\oplus I$ and is defined by the set $\{(r,i): r\in R,i\in I \}$ with addition is component-wise and multiplication is defined by $(r,i)(s,j)=(rs,rj+si+ij)$.
\end{defi}

\begin{lem}
\label{1.6}
\cite[Lemma 4.3]{ggr}
Let $(R,\mathfrak{m})$ be a local ring with maximal ideal $\mathfrak{m}$. Then the excision ring $R\oplus I$ with respect to a proper ideal $I$ in $R$ is also a local ring with maximal ideal $\mathfrak{m}\oplus I$.
\end{lem}

\begin{defi}
We shall say a ring homomorphism $\phi: B\rightarrow D$ is a retract if there exists a ring homomorphism $\gamma: D\hookrightarrow B$ so that $\phi\circ \gamma$ is identity on $D$. We shall also say that $D$ is a retract of $B$.
\end{defi}

\begin{lem}
\label{3.2}
$($\cite[Lemma 3.3]{acr}$)$
Let $B,D$ be rings and let $D$ be a retract of $B$ and let $\pi: B \twoheadrightarrow D$. If $J = \ker(\pi)$, then $E_n(B, J) = E_n(B)\cap SL_n(B, J)$, $n\geq3$.
\end{lem}

\begin{prop}
\label{1.28}
$($\cite[Proposition 3.1]{keshari}$)$
Let $R$ be a commutative ring and $I \subset R$ be an ideal. Then the excision ring $R\oplus I$ is the fibre product of $R$ and $R$ with respect to $R/I$. In fact, if $\phi: D \rightarrow R\oplus I$ is defined by  $\phi(a,b)=(b,a-b)$, then $\phi$ is an isomorphism.
\end{prop}

For all $n\in \mathbb{N}$, let us define a map $i_{2n}: Alt^{\prime}_{2n}(R,I)\ra W^{\prime}_E(D)$ by $i_{2n}(\alpha)=(\chi_r^{\pm 1}\perp \chi_{n-r}^{\pm 1},\alpha)$  where $\ga \equiv \chi^{\pm 1}_r\perp \chi^{\pm 1}_{n-r} \pmod I$ for some $r\geq 0$. It can be shown that the maps $i_{2n}$ will induce a map $i: Alt^{\prime}(R,I)\ra W^{\prime}_E(D)$. It can be shown that $i$ induce a well defined map between $W^{\prime}_E(R,I)$ and $W^{\prime}_E(D)$, which we still call $i$.

\begin{theo}
Let $R$ be a commutative ring and $I\subset R$ be an ideal. Then The set $W^{\prime}_E(R,I)$ has an abelian group structure with the operation $[\alpha].[\beta]=[\alpha\perp\beta]$.
\end{theo}

${\pf}$ It is enough to check that each element $[\alpha]$ of $W^{\prime}_E(R,I)$ has an inverse. Let $\alpha\in Alt^{\prime}_{2n}(R,I)$. Consider the element $\tilde{\alpha}$ defined by $\phi\circ i_{2n}(\alpha)$. We have $\tilde{\alpha}\equiv \chi^{\pm 1}_r\perp \chi_{n-r}^{\pm 1} \pmod{0\oplus I}$ for some $r\geq 0$. In the group $W^{\prime}(R\oplus I)$, we have $[\tilde{\alpha}\perp \tilde{\alpha}^{-1}]=[\chi_1]$. Hence there exists $E\in E_{4n+2t}(R\oplus I)$ such that $E^T(\tilde{\alpha}\perp \tilde{\alpha}^{-1}\perp \chi_t)E=\chi_{2n+t}$.  Going modulo $0\oplus I$ we have $\bar{E}^T(\chi^{\prime} \perp (\chi^{\prime})^{-1}\perp \chi_t)\bar{E}= \chi_{2n+t}$ where $\chi^{\prime}= \chi^{\pm 1}_r\perp \chi^{\pm 1}_{n-r}$. Since $[\chi_r^{\pm 1}\perp \chi_{n-r}^{\pm 1}]=[\chi_n]$ in $W^{\prime}(R,I)$, then there exist $F\in E_{4n+2t}(R,I)$ such that  $F^T(\chi^{\prime}\perp (\chi^{\prime})^{-1}\perp \chi_t) F= \chi_{2n+t}$. Also we have $\tilde{F}^T(\chi^{\prime}\perp (\chi^{\prime})^{-1}\perp 
\chi_t) \tilde{F}= \chi_{2n+t}$ where $\tilde{F}\in E_{4n+2t}(R\oplus I, 0\oplus I)$.  Now replacing $E$ by $E\bar{E}^{-1} \tilde{F}$, we may assume that $E\in E_{4n+2t}(R\oplus I, 0\oplus I)$ and $E^T(\tilde{\alpha}\perp \tilde{\alpha}^{-1}\perp \chi_t)E= \chi_{2n+t}$. Now projecting $R\oplus I$ onto $R$, we have $E_1^T(\alpha\perp \alpha^{-1}\perp \chi_t)E_1= \chi_{2n+t}$ where $E_1\in E_{4n+2t}(R,I)$. Hence we have $[\alpha\perp\alpha^{-1}]=[\chi_1]$, i.e., $[\alpha].[\alpha^{-1}]=[\chi_1]$. Similarly we can show that $[\alpha^{-1}].[\alpha]=[\chi_1]$. 
Hence $W^{\prime}(R,I)$ is a group with the operation $[\alpha].[\beta]=[\alpha\perp \beta]$. To show this group structure is commutative: Let $\alpha\in Alt_{2n}(R,I)$ and $\beta \in Alt_{2m}(R,I)$. Hence we have $\alpha\equiv \chi_r^{\pm 1}\perp\chi_{n-r}^{\pm 1} \pmod I$ and $\beta \equiv \chi_s^{\pm 1}\perp \chi_{m-s}^{\pm 1} \pmod I$ for some $r\geq 0, s\geq 0$.
Now consider the element $\tilde{\alpha}, \tilde{\beta}$ in $Alt_{2n}(R\oplus I)$ and $Alt_{2m}(R\oplus I)$ respectively. In $W^{\prime}(R\oplus I)$, we have $[\tilde{\alpha}\perp \tilde{\beta}\perp \tilde{\alpha}^{-1}\perp \tilde{\beta}^{-1}]= [\chi_1]$. Hence there exist $t\in \mathbb{N}$ and $E\in E_{4(m+n)+2t}(R\oplus I)$ such that $E^T(\tilde{\alpha}\perp \tilde{\beta}\perp \tilde{\alpha}^{-1}\perp\tilde{\beta}^{-1}\perp \chi_t)E=\chi_{2(m+n)+t}$. Now going modulo $0\oplus I$, we have $\bar{E}(\chi^{\prime}\perp \chi^{\prime\prime}\perp (\chi^{\prime})^{-1}\perp (\chi^{\prime\prime})^{-1}\chi_t)\bar{E}=\chi_{2(m+n)+t}$ where $\chi^{\prime}= \chi^{\pm 1}_r\perp\chi^{\pm 1}_{n-r}$ and $\chi^{\prime\prime}= \chi^{\pm 1}_s\perp \chi^{\pm 1}_{m-s}$. Now there exist $F\in E_{4(m+n)+2t}(R\oplus I, 0\oplus I)$ such that $F^T(\chi^\prime\perp \chi^{\prime\prime})\perp (\chi^\prime)^{-1}\perp (\chi^{\prime\prime})^{-1}\perp \chi_t)F= \chi_{2(m+n)+t}$. Hence replacing $E$ by $E\bar{E}^{-1}F$ we may assume that $E\in E_{4(m+n)+2t}(R\oplus I,0\oplus I)$ and $E^T(\tilde{\alpha}\perp \tilde{\beta}\perp \tilde{\alpha}^{-1}\perp \tilde{\beta}^{-1}\perp \chi_t)E= \chi_{2(m+n)+t}$. Now projecting $R\oplus I$ onto $R$, we have $E_1^T(\alpha\perp \beta \perp \alpha^{-1}\perp \beta^{-1}\perp \chi_t)E_1= \chi_{2(m+n)+t}$ where $E_1\in E_{4(m+n)+2t}(R,I)$. Hence $[\alpha\perp \beta \perp \alpha^{-1}\perp \beta^{-1}]=[\chi_1]$ in $W^{\prime}(R,I)$. Hence $[\alpha\perp \beta]=[\beta \perp \alpha]$. Hence $W^{\prime}(R,I)$ is an abelian group with the operation $\perp$.
$~~~~~~~~~~~~~\qedwhite$

\begin{lem}
\label{1.23}
 The following sequence is exact
 \[
\begin{tikzcd}
0 \arrow{r}{} &W^{\prime}_E(R,I) \arrow{r}{i} &W^{\prime}_E(D) \arrow{r}{(p_1)_{*}} &W^{\prime}_E(R) \arrow{r}{} &0
\end{tikzcd}
\]
 Where $i([\alpha])=[(\chi_r^{\pm 1}\perp \chi^{\pm 1}_{n-r},\alpha)]$ and $p_1:D\rightarrow R$ is the projection onto first factor.
\end{lem}

${\pf}$ Let $i([\alpha])=[\chi_1]$ in $W^{\prime}_E(D)$. Consider the isomorphism between $W^{\prime}_E(D)$ and $W^{\prime}_E(R\oplus I)$ induced from the isomorphism $\phi: D\rightarrow R\oplus I$ given by $\phi(a,b)=(b,a-b)$. Then we have $\phi\circ i([\alpha])=[\chi_1]$. Hence  we have $[\tilde{\alpha}]=[\chi_1]$ where $\tilde{\alpha}=\phi((\chi_r^{\pm 1}\perp \chi^{\pm 1}_{n-r},\alpha))$. We have $\tilde{\ga} \equiv \chi_r^{\pm 1}\perp \chi^{\pm 1}_{n-r} \pmod{0\oplus I}$. Replacing $\tilde{\ga}$ by $F^T\tilde{\ga}F$ for some $F\in E_{2n}(R\oplus I, 0\oplus I)$, we may assume that $\tilde{\ga}\equiv \chi_n \pmod{0\oplus I}$. Since $[\tilde{\ga}]=[\chi_1]$ in $W^{\prime}_E(R\oplus I)$, there exists  a natural number $t$ and an elementary matrix $E\in E_{2(n+t)}(R\oplus I)$ such that
 \begin{center}
$E^T(\tilde{\alpha} \perp \chi_t)E=\chi_{n+t}$.
\end{center}
 Now going modulo $0\oplus I$, we have, $\bar{E}^T\chi_{n+t}\bar{E}= \chi_{n+t}$. Hence replacing $E$ by $E\bar{E}^{-1}$, we may assume that $E^T(\tilde{\ga} \perp \chi_t)E= \chi_{n+t}$ and  $E\in E_{2(n+t)}(R\oplus I, 0\oplus I)$ by Lemma \ref{3.2}. Hence projecting $R\oplus I$ onto $R$, we have $\varepsilon^T(\alpha\perp \chi_t)\varepsilon=\chi_{n+t}$, where $\varepsilon\in E_{2(n+t)}(R,I)$. Hence $i$ is injective.

Clearly $(p_1)_*$ is surjective. Now let $(p_1)_*([\beta])=[\chi_1]$. By Lemma \ref{1.5}, we have $\beta=(\beta_1,\beta_2)$ with $\beta_1 \equiv \beta_2 \pmod I$. Thus $[\beta_1]=[\chi_1]$ in $W^{\prime}_E(R)$. Hence there exists a natural number and $F\in E_{2(n+t)}(R)$ such that 
\begin{center}
$F^T(\beta_1\perp\chi_t)F=\chi_{n+t}$
\end{center}
Now $(F,F)^T((\beta_1,\beta_2)\perp (\chi_t,\chi_t))(F,F)=(\chi_{n+t},F^T(\beta_2\perp \chi_t)F)$.

Hence $[\beta]=[(\beta_1,\beta_2)]=[(\beta_1,\beta_2)\perp(\chi_t,\chi_t)]=[(\chi_{n+t},F^T(\beta_2\perp\chi_t)F]=i([F^T(\beta\perp \chi_t)F])$.

Hence the sequence
 \[
\begin{tikzcd}
0 \arrow{r}{} &W^{\prime}_E(R,I) \arrow{r}{i} &W^{\prime}_E(D) \arrow{r}{(p_1)_*} &W^{\prime}_E(R) \arrow{r}{} &0
\end{tikzcd}
\]
is exact.
\\
$~~~~~~~~~~~~~~~~~~~~~~~~~~~~~~~~~~~~~~~~~~~~~
\qedwhite$

\begin{lem}
\label{1.22}
\cite{fr}
Let $C$ be the kernel of the group homomorphism $R^{\times}\rightarrow (R/I)^{\times}$ induced by $\pi$, with the convension that $C=R^{\times}$if $I=R$. Then the Pfaffian gives a split exact sequence
 \[
\begin{tikzcd}
0 \arrow{r}{} &W_E(R,I) \arrow{r}{j} &W^{\prime}_E(R,I) \arrow{r}{Pf} &C \arrow{r}{} &0
\end{tikzcd}
\]
\end{lem}
${\pf}$  Clearly the homomorphism $i$ is injective and the sequence is exact on the middle. For split exact sequence we have to show that there is a map $r:C\rightarrow W^{\prime}_E(R,I)$ such that $Pf\circ r=id_C$. Define
\begin{eqnarray*}
r(a)=\left[\begin{pmatrix}
0 &a\\
-a &0
\end{pmatrix}\right]
\end{eqnarray*}
The map $r$ is well-defined since $a\equiv 1\pmod I $.
\\
$~~~~~~~~~~~~~~~~~~~~~~~~~~~~~~~~~~~~
\qedwhite$

\begin{cor}
The set $W_E(R,I)$ is a group with respect to the operation $\perp
$.
\end{cor}
\begin{lem}
\label{1.24}
 The following sequence is exact.
 \[
\begin{tikzcd}
0 \arrow{r}{} &W_E(R,I) \arrow{r}{i} &W_E(D) \arrow{r}{(p_1)_*} &W_E(R) \arrow{r}{} &0
\end{tikzcd}
\]
\end{lem}

${\pf}$ Consider the following diagram
\[
\begin{tikzcd}
&0 \arrow{d}{}
               &0 \arrow{d}{} &0 \arrow{d}{}\\
0\arrow{r}{} &W_E(R,I) \arrow{r}{j} \arrow{d}{i} &W^{\prime}_E(R,I) \arrow{r}{Pf} \arrow{d}{i}&C\arrow{r}{}\arrow{d}{\phi} &0\\
0 \arrow{r}{} &W_E(D)\arrow{r}{j}\arrow{d}{(p_1)_*} &W^{\prime}_E(D)\arrow{r}{Pf}\arrow{d}{(p_1)_*} &D^{\times}\arrow{r}{} \arrow{d}{\psi}&0\\
0 \arrow{r}{} &W_E(R) \arrow{r}{j}\arrow{d}{} &W^{\prime}_E(R) \arrow{r}{Pf}\arrow{d}{}&R^{\times} \arrow{r}{}\arrow{d}{}&0\\
&0 &0 &0
\end{tikzcd}
\]
Where the maps $\phi:C\rightarrow D^{\times}$ and $ \psi: D^{\times}\rightarrow R^{\times}$ are defined by $\phi(a)=(1,a)$ and $\psi(d_1,d_2)=d_1$.
With respect to the $\phi$ and $\psi$ the above diagram is commutative. Each row of the above diagram is an exact sequence by Lemma \ref{1.22} . Observe that all these groups in the above diagram are abelian. The  second column of the diagram is an exact sequence by Lemma \ref{1.23}. Also the third column of the above diagram is exact. Hence by diagram chasing the first column is also exact.
$~~~~~~~~~~~~~~~~~~~~~~~~~~~~~~~~~~~~~~~~~~~~~~
\qedwhite$

\section{\bf{Divisibility of $W_E(R[X],I[X])$}}

\begin{lem}
\label{1.14}
$($\cite[Lemma 3.1]{SV}$)$ $($Karoubi$)$ Let $R$ be a commutative ring with $1$. Let $\alpha \in W_E(R[X])$. Then we have $[\alpha]=[l]$, where $l=\varphi_0+\varphi_1X$ where $\varphi_0$ and $\varphi_1$ are matrices over $R$.
\end{lem}

\begin{lem}
\label{1.14.0}
$($Karoubi$)$ Let $R$ be a commutative ring and $I\subset R$ be an ideal of $R$. Let $\alpha\in W_E(R[X],I[X])$. Then we have $[\alpha]=[l]$, where $l=\varphi_0+\varphi_1X$  where $\varphi_0$ is a matrix over $R$ and $\varphi_1$ is a matrix over $I$.
\end{lem}

${\pf}$
Let $\alpha\in Alt_{2n}(R[X],I[X])$. We may assume that $\alpha\equiv \chi_n \pmod{I[X]}$. Consider the element $\tilde{\alpha}\in Alt_{2n}((R\oplus I)[X])$. Then we have $\tilde{\alpha}\equiv \chi_n \pmod{0\oplus I[X]}$. Now by Karoubi we have $[\tilde{\alpha}]=[l]$, where $l$ is a linear matrix in $Alt((R\oplus I)[X])$. It is easy to check that $l\in Alt_{2m}((R\oplus I)[X],0\oplus I[X])$ for some $m$. We have $[\tilde{\alpha}\perp l^{-1}]=[\chi_1]$. Conjugating $\tilde{\alpha}\perp l^{-1}$ with some element $F\in E_{2(m+n)}((R\oplus I)[X], 0\oplus I[X])$, we may assume that $\tilde{\alpha}\perp l^{-1}\equiv \chi_{m+n} \pmod{0\oplus I[X]}$. Since $[\tilde{\alpha}\perp l^{-1}]=[\chi_1]$, then there exist $t\in \mathbb{N}$ and $E\in E_{2(m+n+t)}((R\oplus I)[X])$ such that $E^T(\tilde{\alpha}\perp l^{-1}\perp \chi_t)E=\chi_{m+n+t}$. Now going modulo $0\oplus I[X]$, we have $\bar{E}^T\chi_{m+n+t}\bar{E}=\chi_{m+n+t}$. Hence replacing $E$ by $E\bar{E}^{-1}$, we have $E\in E_{2(m+n+t)}((R\oplus I)[X],0\oplus I[X])$ and $E^T(\tilde{\alpha}\perp l^{-1}\perp \chi_t)E= \chi_{m+n+t}$. Now projecting $(R\oplus I)[X]$ onto $R[X]$, we have $E_1^T(\alpha \perp l_1^{-1} \perp \chi_t)E_1=\chi_{m+n+t}$, where $E_1\in E_{2(m+n+t)}(R[X],I[X])$ and $l_1\in Alt_{2m}(R[X],I[X])$ is a linear matrix. Hence we have $[\alpha\perp l_1^{-1}]=[\chi_1]$ in $W_E(R[X], I[X])$. In other words $[\alpha]=[l_1]$ in $W_E(R[X],I[X])$. 
$~~~~~\qedwhite$

\begin{prop}
\label{1.15}
\cite[Proposition 2.4.1]{rao}
Let $R$ be a local ring with $1/2k\in R$, and let $[V]\in W_E(R[X])$. Then $[V]$ has a $k$-th root, i.e. there is a $[W]\in W_E(R[X])$ such that $[V]= [W]^k$ in $W_E(R[X])$.
\end{prop}

\begin{lem}
Let $R$ be a local ring and $I\subset R$ be an ideal of $R$, $1/2k\in R$ and let $[\alpha]\in W_E(R[X],I[X])$.  Then $[\alpha]$ has a $k$-th root.
\end{lem}

${\pf}$ Let $I=R$. By convension we have $W_E(R[X], R[X])= W_E(R[X])$. The lemma is true for $I=R$ by Proposition \ref{1.15}. Hence we may assume that $I$ is a proper ideal of $R$. Let $\alpha \in Alt_{2n}(R[X], I[X])$.  We may assume that $\alpha\equiv \chi_n (\mod I)$. Consider the element $\tilde{\alpha}\in W_E((R\oplus I)[X])$. By Lemma \ref{1.6}, we have $R\oplus I$ is a local ring. And also we have $1/2k\in R\oplus I$. Hence by Proposition \ref{1.15}, we have there exists $\beta\in W_E((R\oplus I)[X])$ such that $[\tilde{\alpha}]=[\beta]^k$. We revisit the proof of Proposition \ref{1.15}, to verify  that $\beta\equiv \chi_m \pmod{0\oplus I[X]}$ for some $m$. We have $\tilde{\alpha}\in Alt_{2n}((R\oplus I)[X],0\oplus I[X])$.  
 By Lemma \ref{1.14.0}, we have there exist  $t\in \mathbb{N}$ and $E\in E_{2(n+t)}((R\oplus I)[X], 0\oplus I[X])$ such that $E^T(\tilde{\alpha}\perp \chi_t)E= \chi_{n+t}+rX$, where $r$
 is a matrix over $0\oplus I$. Now let $\gamma=I_{2(n+t)}-\chi_{n+t}rX$. Clearly $\gamma\in SL_{2(n+t)}((R\oplus I)[X],0\oplus I[X])$. Hence $\chi_{n+t}r$ is a nilpotent matrix over $R\oplus I$. Since $1/2k\in R\oplus I$, we can extract a $2k$-th root of $\gamma$. Call it $\delta$. It is easy to check that $\delta\in SL_{2(n+t)}((R\oplus I)[X], 0\oplus I[X])$. Now M. Karoubi pointed out that $$E^T(\tilde{\alpha}\perp \chi_t)E= \chi_{n+t}\gamma= \chi_{n+t}\delta^{2k}=(\delta^k)^T\chi_{n+t}\delta^k.$$
Let $\beta= \delta^T\chi_{n+t}\delta$. Clearly $\beta\equiv \chi_{n+t} (\mod 0\oplus I[X])$. Applying Whitehead's lemma one can check that $[\tilde{\alpha}]=[\beta]^k$. Let $m= n+t$.  We have \begin{eqnarray}
[\tilde{\alpha}\perp \beta^{-1}\perp .. ~\text{(k times)}~ .. \perp \beta^{-1}]=[\chi_1]
\end{eqnarray}
 in $W_E((R\oplus I)[X])$. Conjugating $\tilde{\alpha}\perp \beta^{-1}\perp \dots \perp \beta^{-1}$ by some element in $E_{2(n+mk)}((R\oplus I)[X],0\oplus I[X])$ we may assume that $\tilde{\alpha}\perp \beta^{-1}\perp \dots \perp \beta^{-1} \equiv \chi_{n+mk} \pmod{0\oplus I[X]}$. By $(1)$ we have there exist $p\in \mathbb{N}$ and $F\in E_{2(n+mk+p)}((R\oplus I)[X])$ such that $F^T(\tilde{\alpha}\perp \beta^{-1}\perp \dots \perp \beta^{-1}\perp \chi_p)F=\chi_{n+mk+p}$. Now going modulo $0\oplus I[X]$, we have $\bar{F}^T\chi_{n+mk+p}\bar{F}=\chi_{n+mk+p}$. Hence replacing $F$ by $F\bar{F}^{-1}$, we have $F\in E_{2(n+mk+p)}((R\oplus I)[X],0\oplus I[X])$ and $F^T(\tilde{\alpha}\perp \beta^{-1}\perp \dots \perp \beta^{-1}\perp \chi_p)F=\chi_{n+mk+p}$. Now projecting $R\oplus I$ onto $R$ we have $F_1^T(\alpha\perp \beta_1^{-1}\perp \dots \perp \beta_1^{-1}\perp \chi_p)F_1=\chi_{n+mk+p}$, where $F_1\in E_{2(n+mk+p)}(R[X],I[X])$. Hence we have $[\alpha\perp \beta_1^{-1}\perp \dots \perp \beta_1^{-1}]=[\chi_1]$ in $W_E(R[X],I[X])$. Hence $[\alpha]= [\beta_1]^k$ in $W_E(R[X],I[X])$.
 $~~~~~~~~~~~\qedwhite$
 
\begin{lem}
Let $R$ be a local ring of dimension $d$, $d>1$, $\alpha\in R$ be a non-zero-divisor, $v\in Um_{d+1}(R[X], (\alpha))$. Then $v\sim_{E_{d+1}(R[X],(\alpha))} (a_0(X),\alpha a_1(X), \dots, \alpha c_d)$ such that $c_d$ is a non-zero-divisor.
\end{lem}

${\pf}$ By Excision theorem \cite[Theorem 3.21]{vdk1}, we may assume that, $R$ is a reduced ring. Let $v(X):=(v_0(X),v_1(X),\dots, v_d(X))$. We assume that $\deg v_0(X)\geq 1$. Let the leading coefficient of $v_0(X)$ be $a$. We may assume that $a$ is a non-zero-divisor of $R$. Let the 'overline'  denote modulo $(a)$ and consider $\bar{v(X)}\in Um_{d+1}(\bar{R}[X],(\alpha))$. By excision and usual stability estimates we have $\bar{v(X)}\sim_{E_{d+1}(\bar{R}[X],(\alpha))} e_1$. Hence, we may modify $v(X)$ suitably and assume that, $\bar{v(X)}=e_1$. As, so crucially, indicated by M. Roitman, this transformation can be performed so that every stage the row contains a polynomial which is unitary in $R_a$. And we may ensure that $v_0(X)$ is unitary in $R_a$ and  $\deg v_0(X)\geq 1$. Let $v_0(X)=1+v_0^{\prime}(X)$, $v_i(X)=\alpha av_i^{\prime}(X)$ for $i>0$. One has
 \begin{center}
$ a^{l_i}v_i^{\prime}(X)= q_i(X)v_0(X)+r_i(X)$,
 \end{center} 
 for some $l_i>0$, $r_i(X)\in R[X]$, with $r_i(X)=0$ or  $\deg r_i(X)< \deg v_0(X)$.  By Excision theorem \cite[Theorem 3.21]{vdk1}, we can transform $v(X)$ in the relative orbit with respect to $(a\alpha)$ and assume that $\deg v_i^{\prime}(X)< \deg v_0(X)$, for all $i>0$. If $\deg v_0(X)=1$, then we are done. Hence we may assume that $\deg v_0(X)=d_0\geq 2$. 
 
 Let $c_1,c_2,\dots, c_{d_0(d-1)}$ be the coefficients of $1,X,\dots, X^{d_0-1}$ of the polynomials $v_2^{\prime}(X),\dots, v_d^{\prime}(X)$. Since $d_0(d-1)\geq 2.\frac{(d+1)}{2}> \dim R_a$, we can argue as in the proof of \cite[Theorem 5]{roitman} to conclude that  the ideal generated by the polynomials $v_0(X),v_2^{\prime}(X), \dots, v_d^{\prime}(X)$ contains a polynomial $h(X)$ of degree $(d_0-1)$  which is unitary in $R_a$. Let  leading coefficient of $h(X)$ is $ua^k$ where $u$ is a unit in $R$. Via Excision theorem and the argument in the proof of \cite[Theorem 5]{roitman}, we have 
 
\begin{align*}
 & mse(v_0(X),\alpha v_1^{\prime}(X),\dots, \alpha v_d^{\prime}(X))\\
 &=mse(v_0(X), a^{2k}\alpha v_1^{\prime}(X),\dots, \alpha v_d^{\prime}(X)) \\
 &= mse(v_0(X),\alpha \{a^{2k}v_1^{\prime}(X)+(1-a^k.l(v_1^{\prime}(X)))h\}, \alpha v_2^{\prime}(X),\dots, \alpha v_d^{\prime}(X))\\
 &=mse(v_0(X),\alpha v^{\prime\prime}(X), \alpha v_2^{\prime}(X),\dots, \alpha v_d^{\prime}(X))
 \end{align*}
 where $l(v_1^{\prime}(X))$ is the leading coefficient of $v_1^{\prime}(X)$. Note that $v_1^{\prime\prime}(X)$ is  unitary in $R_a$ of degree $d_1<d_0$. Rename $v_1^{\prime\prime}(X)$ as $v_1^{\prime}(X)$. Again we may ensure that $\deg v_i^{\prime}(X) <d_1$, for $i\geq 2$. Repeat the argument above to lower the degree of $v_1^{\prime}(X)$ to $1$. Then lower the degree of $v_i^{\prime}$ to zero. We then have the desired form of the vector $v(X)$ in the class of $MSE_{d+1}(R[X],(\alpha))$. 
 
 $~~~~~~~~~~~~~~~~~~~~~~~~~~~~~~~~~~~~~~~~~~~~~~~~\qedwhite$
 
 \section{\bf{Analytic isomorphisms}}
 \begin{defi}
Let $A, B$ be  commutative rings  and $\phi: B\rightarrow A$ be a ring homomorphism. Let $h\in B$. $\phi$ is said to be analytic  isomorphism along $h$, if the following conditions are satisfied.

(i) $h$ is a non-zero-divisor of $B$.

(ii) $\phi(h)$ is a non-zero-divisor of $A$.

(iii) $\phi$ induces an isomorphism between $B/hB$ and $A/hA$.
\end{defi}

In this situation the commutative diagram
 \begin{equation}\label{eqn:0}
\begin{tikzcd}
B \arrow{r}{\phi}\arrow{d}{} &A \arrow{d}{}\\
B_h\arrow{r}{\phi_h} &A_h
\end{tikzcd}
\end{equation}

will be called an analytic diagram. In \cite{roy}, it is checked that $(1)$ is a cartesian square. 

Let $\textbf{P}(R)$ denote the category of all finitely generated projective $R$-modules. In \cite{roy}, it is shown that the corresponding square for projective modules is also cartesian, i.e.,

\[
\begin{tikzcd}
\textbf{P}(B)\arrow{r}{\phi}\arrow{d}{}&\textbf{P}(A)\arrow{d}{}\\
\textbf{P}(B_h)\arrow{r}{\phi_h}&\textbf{P}(A_h)
\end{tikzcd}
\]
 
is cartesian.

\begin{lem}
\label{1.38}
Let $A$ be a commutative ring and $B\subset A$ be a subring. Let $h\in B$ be such that $B\subset A$ is analytically isomorphic along $h$. Let $I$ be an ideal of $A$. Then $B\oplus I\subset A\oplus I$ is analytically isomorphic along $(h,0)$.
\end{lem}
${\pf}$
Claim: $(h,0)$ is a non-zero-divisor of $B\oplus I$.

Let $(h,0)(h_1,i_1)=(0,0)$ in $B\oplus I$ for some $(h_1,i_1)\in B\oplus I$. Then $(hh_1,hi_1)=(0,0)$. Thus we have $hh_1=0$, $hi_1=0$. Since $h$ is a non-zero-divisor of $B$ and of $A$ also, then $h_1=0$ and $i_1=0$. Hence $(h_1,i_1)=(0,0)$ and therefore $(h,0)$ is a non-zero-divisor of $B\oplus I$.

Similarly  $(h,0)$ is a non-zero-divisor of $A\oplus I$. 

Consider the natural map $\phi:(B\oplus I)/(h,0)(B\oplus I)\rightarrow (A\oplus I)/(h,0)(A\oplus I)$. Let $\phi(\overline{(b,i)})=0$. Then $(b,i)\in (h,0)(A\oplus I)$. Hence $(b,i)=(h,0)(a,j)$ for some $(a,j)\in A\oplus I$. This gives $b=ha$ and $i=hj$. Since $B\subset A$ is analytically isomorphic along $h$, $b=ha$ implies that there exists $b^{\prime}\in B$ such that $b=hb^{\prime}$. Hence we have $(b,i)=(h,0)(b^{\prime},j)$. This shows that $\phi$ is injective. Let $\overline{(a,i)}\in (A\oplus I)/(h,0)(A\oplus I)$. Since $B\subset A$ is analytic along $h$, then there exists $b\in B$ and $c\in A$ such that $a-b=hc$. This shows that $\phi(\overline{(b,i)})=\overline{(a-hc,i)}= \overline{(a,i)}$. Hence $\phi$ is surjective. Thus $B\oplus I\subset A\oplus I$ is analytically isomorphic along $(h,0)$.
$~~~~~~~~~~~~~~~~~~~~~~~~~~~~~~~~~~~~~~~~~~~~~~~~
\qedwhite$

\begin{defi}
Let $(R,\mathfrak{m})$ be a local ring. A monic polynomial $f\in R[X]$ is said to be Weierstrass polynomial if $f= X^n+ a_1X^{n-1}+...+a_n, a_i \in \mathfrak{m}~ for~ i=1,2,...,n$.
\end{defi}

\begin{lem}
\label{0.6}
$($\cite[Proposition 1.7]{budh}$)$
Let $(R,\mathfrak{m})$ be a local ring and $f\in R[X]$ be a   Weierstrass polynomial. Then $R[X] \hookrightarrow R[X]_{(\mathfrak{m},X)}$ is an analytic isomorphism along $f$.
\end{lem}

\begin{lem}
\label{1.37}
$($\cite [Lemma 2.4] {Vorst}$)$
Let $B \subset A$ be a subring, $h\in B$ be such that $B\hookrightarrow A$ is an analytic isomorphism along $h$. Let $\alpha \in E_r(A_h)$, $r\geq 3$. Then there exists $\alpha_1 \in E_r(A)$ and $\alpha_2 \in E_r(B_h)$ such that $\alpha =(\alpha_1)_h (\alpha_2\otimes 1)$.

 Consequently if $\alpha\in GL_r(A)$ with $\alpha_h\in E_r(A_h)$ then there exist  $\alpha_1 \in E_r(A)$  and  $\alpha_2 \in GL_r(B)$  such that $(\alpha_2)_h\in E_r(B_h)$ and $\alpha=(\alpha_1)_h (\alpha_2\otimes 1)$.
\end{lem}

\begin{defi}
Let $R$ be a local ring. $R$ is said to be a local algebra with a ground field if $R$ is a localisation of an affine $k$-algebra for some field $k$.
\end{defi}

\begin{defi}
Let $R$ be a regular local algebra with a ground field. $R$ is said to be regular local algebra with a separating ground field if $R$ possesses a ground field $K$ such that the residue field of $R$ is a finite separable and hence a simple extension of $K$.
\end{defi}

\begin{defi}
\textbf{Regular $k$-spots}
 Let $k$ be a field. By a regular spot over a field  $k$ we mean a localisation $C_{\mathfrak{p}}$ of a finitely generated $k$-algebra C at a regular prime $\mathfrak{p} \in  \Spec C$.
\end{defi}

Lindel \cite[Proposition 2]{lindel} analysed regular k-spots over a perfect field as \'{e}tale extensions of rings of the type $k[X_1,...,X_n]_{(\varphi(X_1),X_2,...,X_n)}$, where $\varphi(X_1)\in k[X_1]$ is a monic polynomial. %And we can resolve the case for regular $k$-spot over arbitrary field by argument given in \cite{lindel}  by Swan. 

\begin{lem}
\label{1.35}
\cite[Proposition 6.10]{ak}
Let $A$ be a noetherian ring, $B$ an $A$-algebra, $\mathfrak{q}$ a prime ideal of $B$, $\mathfrak{p}$ the trace of $\mathfrak{q}$ in $A$. Suppose there exists a polynomial $P(T)$ and an element $t\in B$ such that the map $A[T]/PA[T]\rightarrow B$ defined by $t$ is an isomorphism. Then $B_{\mathfrak{q}}$ is unramified over $A_{\mathfrak{p}}$ if and only if $(P(T),P^{\prime}(T))A[T]=A[T]$. Suppose, in addition, that the leading coefficient of $P$ is invertible. Then $B_{\mathfrak{q}}$ is \'{e}tale over $A_{\mathfrak{p}}$ if and only if $P^{\prime}(t)\notin \mathfrak{q}$.
\end{lem}

\begin{theo}
\label{1.36}
$($\cite[Theorem 2.8]{budh}$)$
Let $(R,\mathfrak{m})$ be a regular $k$-spot of dimension $d$ with a separating ground field $K$. Let $g$ be any element of $\mathfrak{m}^2$.
Then there exists a regular local subring $S$ of $R$ such that:

(i) $S = K[X_1,..., X_d]_{(\varphi(X_1),X_2,\dots,X_d))}$, where $\varphi(X_1) \in K[X_1]$ is an irreducible monic polynomial.

(ii) $S\subset R$ is an analytic isomorphism along $h$, for some $h \in gR\cap S$.(Here $h$ depends on choice of $g$.)
 
\end{theo}

We revisit the proof of \cite[Theorem 2.8]{budh} to justify the following crucial (in our context) claim:
\begin{lem}
\label{1.37}
The element $h$ and $g$ mentioned in Theorem \ref{1.36} differ by a unit in $R$.
\end{lem}

${\pf}$ The regular local subring $S$ of $R$ is constructed in the Theorem \ref{1.36} in the following way:

We can choose elements $X_2,X_3,\dots,X_d$ in $\mathfrak{m}$ such that $g,X_2,X_3,\dots,X_d$ is a system of generators in $R$ and $X_2,X_3,\dots,X_d$ are a part of a minimal generating set for $\mathfrak{m} \pmod{\mathfrak{m}^2}$. Since $R$ is a regular local ring we have that $g,X_2,\dots X_d$ is a regular sequence in $R$. The field $K$ is contained in $R$. Therefore $a,X_2,\dots, X_d$ are algebraically independent over $K$. Thus $C_1= K[g,X_2,\dots,X_d]$ is a polynomial ring contained in $R$. Let $B$ be the integral closure of $C_1$ in $R$ and let $\mathfrak{m}_1=\mathfrak{m}\cap B$. Since $\mathfrak{m}_1$ contracts to the maximal ideal $(g,X_2,\dots,X_d)$ of $C_1$, it follows that $\mathfrak{m}_1$ is a maximal ideal in $B$. It is proved in \cite[Theorem 2.8]{budh} that $R= B_{\mathfrak{m}_1}$. Viewing $R$ as $B_{\mathfrak{m}_1}$ we may rename $\mathfrak{m}_1$ as $\mathfrak{m}$. 

By the hypothesis of Theorem \ref{1.36}, we have, $L$:= $Q(R)$= $B/\mathfrak{m}$= $K(\bar{\alpha})$ for some $\alpha$ in $B$. Let $\varphi$ denote the minimal polynomial of $\bar{\alpha}$ over $K$. Then $\varphi(\alpha)\in \mathfrak{m}$. We can choose $X_1=\alpha$ in such a way that $(g,\varphi(X_1),X_2,\dots,X_d)=\mathfrak{m}$.

Further, as $X_1$ is integral over $K[g,X_2,\dots,X_d]$, replacing $X_1$ by $X_1+g^j$ for some integers $j$ large enough, we may even assume that $g$ is integral over $K[X_1,X_2,\dots,X_d]$. We now define $C=K[X_1,X_2,\dots,X_d]$ and $M=(\varphi(X_1),X_2,\dots,X_d)$. Set $S= C_M$.

The element $h$ is chosen in \cite[Theorem 2.8]{budh} in the following way:
 
  Note that $\mathfrak{m}\cap C=M$. Consider the ring $C[g]$ and the maximal ideal $(M,g)$ of $C[g]$. As $(M,g)$ generates $\mathfrak{m}$ in $B$, $B$ is a finite $C[g]$-module and $B/\mathfrak{m}=L=C[g]/(M,g)$, we conclude, using Nakayama's lemma, that $R= B_{\mathfrak{m}}= C[g]_{(M,g)}$. It is proved in \cite[Theorem 2.8]{budh}, that $R$ is faithfully flat $S$-algebra and also $R$ is unramified over $S$. Therefore $R$ is \'{e}tale over $S$. Consider the $C$-algebra homomorphism $\sigma: C[T]\rightarrow C[g]$ by $\sigma(T)=g$. As $g$ is integral over $C$, there exists an irreducible monic polynomial in $C[T]$: $$F(T)= T^n+\lambda_{n-1}T^{n-1}+\dots+\lambda_1T+\lambda_0$$ such that $\sigma(F(T))= F(g) =0$. Thus $\sigma$ induces an isomorphism between $C[T]/FC[T]$ and $C[g]$. Since $R$ is \'{e}tale over $S$, then by Lemma \ref{1.35}, we have $F^{\prime}(g)\notin (M,g)$. Hence we have $\lambda_1\notin M$. Thus $\lambda_0= -g(\lambda_1+\lambda_2g+\dots+\lambda_{n-1}g^{n-1}+g^{n-1})$. In other words $\lambda_0R= gR$ and $\lambda_0\in S\cap gR$. Now take $h=\lambda_0$. Since $R/gR= C[g]_{(M,g)}/gC[g]_{(M,g)}= S/S\cap gR$ we have $R=S+gR= S+hR$. Since $R$ is faithfully flat $S$-algebra we have, $S\cap hR=hS$. Hence we have $S\subset R$ is an analytic isomorphism along $h$.
$~~~~~~~~~~~~~~~~~~~~~~~~~~~~~~~~~~~~~~~~~~~~~~~~~~~~~~~\qedwhite$

\begin{rem}
The following lemma gives a sufficient condition for a local domain to be a local algebra with a separating ground field.
\end{rem}

\begin{lem}
\label{1.39}
$($\cite[Lemma 2.10]{budh}$)$
Let $R$ be a local domain. Suppose that $R$ is a local algebra with a ground field $k$ and that $k$ is a finitely generated field extension of a perfect field $k_0$. Then $R$ is a local algebra with a separating ground field.
\end{lem} 

\begin{lem}
\label{3.5}
$($Nagata$)$
Let $k$ be a field and $f$ be a polynomial in $k[X_1,X_2,...,X_d]$ and $\phi(X_1)\in k[X_1]$ be a monic polynomial. Then there exists a change of variables, $X_1\mapsto X_1$, $X_i\mapsto X_i+\phi(X_1)^{r_i}$, $2\leq i\leq d$, such that $f=c.h(X_1,X_2,...,X_d)$, where $c\in K^*$ and $h$ is a monic polynomial in $X_1$ over $k[X_2,X_3,...,X_d]$.
\end{lem}

${\pf}$ Let $f(X_1,X_2,\cdots,X_d)= \sum_i a_iX_1^{i_1}X_2^{i_2}\cdots X_d^{i_d}$ and $\phi(X_1)=X_1^n+c_1X_1^{n-1}+\cdots+c_n$.
So,\begin{align*}
f(X_1,X_2+\phi(X_1)^{r_2},\cdots,X_d+\phi(X_1)^{r_d}) &=\sum_i a_iX_1^{i_1}(X_2+\phi(X_1)^{r_2})^{i_2}\cdots (X_d+\phi(X_1)^{r_d})^{i_d}\\
&=\sum_i a_iX_1^{i_1+nr_2i_2+\cdots +nr_di_d}+\text{terms of lower degree in } X_1
\end{align*}
Let $m> max\{n.i_k:k=1,\cdots,d\}$. If we choose $r_j=m^{j-1}$, then the integers $i_1+nr_2i_2+\cdots+nr_di_d$ will have different $m-$adic expansion for all tuple $(i_1,i_2,\cdots,i_d)$ in the expansion of $f$. Hence if we choose such $r_i$, then the monomials $a_iX_1^{i_1+nr_2i_2+\cdots+nr_di_d}$ in $f(X_1,X_2+\phi(X_1)^{r_2},...,X_d+\phi(X_1)^{r_d})$ will not cancel out each other. Hence choosing the highest degree with non-zero coefficient among $a_iX_1^{i_1+nr_2i_2+\cdots+nr_di_d}$ we can make it the leading term of $f(X_1,X_2+\phi(X_1)^{r_2},\cdots,X_d+\phi(X_1)^{r_d})$. Hence we have $f=c.h(X_1,X_2,\cdots,X_d)$, where $c\in k^*$ and $h$ is a monic polynomial in $X_1$ over $k[X_2,X_3,\cdots,X_d]$.
$~~~~~~~~~~~~~~~~~~~~~~~~~~~~~~~~~~~~~~~~~~~~~~~~~\qedwhite$

The following is well-known. For completeness we sketch a proof.

\begin{prop}
\label{3.10}
$($\cite[Theorem 2.2.12]{bh}$)$
Let $R^\prime\subset R$ be a faithfully flat local extension. Then $R^\prime$ is regular if $R$ is regular.
\end{prop}

${\pf}$
It is enough to show that $pd_{R^\prime}(M)$ is finite for any finitely generated $R^\prime$-module $M$. Let 
\[
0\rightarrow Q\rightarrow P_m\rightarrow \dots\rightarrow P_0\rightarrow M\rightarrow 0
\] be a resolution of $M$ where $P_i$'s are projective $R^\prime$ modules. Since $R$ is regular we have $pd_R(M\otimes_{R^\prime} R)$ is finite. Hence for large $m$ we have $Q\otimes_{R^\prime} R$ is projective $R$ module. Hence $Q$ is a projective $R^\prime$ module since $R$ is faithfully flat over  $R^\prime$. Hence $pd_{R^\prime}(M)$ is finite.
$~~~~~~~~~~~~~~~~~~~\qedwhite$

\begin{theo}
$($\cite[Corollary 5.7]{S}$)$
\label{3.4}
Let $R$ be a ring and $f\in R[X]$ be a monic polynomial and $\alpha\in GL_r(R[X])$, $r\geq 3$. If $\alpha_f\in E_r(R[X]_f)$, then $\alpha\in E_r(R[X])$.
\end{theo}

\begin{theo}
\label{3.6}
\cite[Theorem 3.12]{kcr}
Let $R$ be a regular $k$-spot of dimension $d$ over a field $k$. Let 
$$\alpha \in SL_n(R[X],(X^2-X))\cap E_{n+1}(R[X],(X^2-X))$$ for  $n\geq 3$, then $\alpha\in E_n(R[X],(X^2-X))$.
\end{theo}

\section{\bf{Relative Vaserstein's Theorem}}
\begin{defi}
A row $(a_1,a_2,...,a_n)\in R^n$ is said to be unimodular if there exist a row $(b_1,b_2,...,b_n)\in R^n$ such that $\sum_{k=1}^{n}a_kb_k=1$. A row $(a_1,a_2,...,a_n)\in R^n$ is said to be relative unimodular row with respect to the ideal $I$ if it is unimodular and $(a_1,a_2,...,a_n)\equiv (1,0,...,0) \pmod I$, i.e. $a_1-1,a_2,...,a_n$  all belong to $I$. The set of all unimodular rows is denoted by $Um_n(R)$ and the set of all relative unimodular rows with respect to the ideal $I$ is denoted by $Um_n(R,I)$.
\end{defi}

\begin{defi}
(Stable range condition $Sr_n(I)$) Let $I$ be an ideal in $R$. We shall say stable range condition $Sr_n(I)$ holds for $I$ if for any $(a_1, a_2, . . . , a_{n+1})\in Um_{n+1}(R, I)$ there exists $c_i $ in $I$ such that $(a_1 +c_1a_{n+1},a_2 +c_2a_{n+1},...,a_n +c_na_{n+1})\in Um_n(R,I)$.
\end{defi}

\begin{defi}
(Stable range $Sr(I)$) We shall define the  stable range of $I$, denoted by $Sr(I)$, to be the least integer $n$ such that $Sr_n(I)$ holds for $I$. When $I=R$, then the stable range of $R$ will be denoted by $Sr(R)$.
\end{defi}

\begin{defi}
(Stable dimension $Sd(R)$) We shall define the stable dimension of R, denoted by $Sd(R)$, to be one less than the stable range of $R$, i.e. $Sd(R)=Sr(R)-1$.
\end{defi}

The group $E_n(R)$ acts on the set $Um_n(R)$ by the action $E_n(R)\times Um_n(R)\longrightarrow Um_n(R)$ given by $\alpha\cdot v=v\alpha$. We denote the orbit space $Um_n(R)/E_n(R)$ by $MSE_n(R)$. Similarly the group $E_n(R,I)$ acts on the set $Um_n(R,I)$ by the same action. We denote the orbit space $Um_n(R,I)/E_n(R,I)$ by $MSE_n(R,I)$.

\subsection*{The Relative Vaserstein symbol}
Let $b=(b_1,b_2,b_3)\in Um_3(R,I)$ and $a=(a_1,a_2,a_3)\in Um_3(R,I)$ be such that $a_1b_1+a_2b_2+a_3b_3=1$.

Denote $\theta(a,b)$ by the matrix
\begin{eqnarray*}
 \theta(a,b)= \begin{pmatrix}
0 & -b_{1} &-b_2 &-b_3\\
b_1&0&-a_3&a_2\\
b_2&a_3&0&-a_1\\
b_3&-a_2&a_1&0
\end{pmatrix} 
\end{eqnarray*}

Define a map $s_E :Um_3(R,I) \rightarrow W_E(R,I)$ by
 $s_G(b)=[\theta(a,b)]$.

\begin{lem}
The map $s_E$ does not depend on the choice of $a$, i.e. if $c=(c_1,c_2,c_3)\in Um_3(R,I)$ be such that $cb^T=ab^T=1$, then $[\theta(c,b)]=[\theta(a,b)]$ in $W_E(R,I)$.
\end{lem}
${\pf}$
Let $\varepsilon=\begin{pmatrix}
1&d_1&d_2&d_3\\
0&1&0&0\\
0&0&1&0\\
0&0&0&1
\end{pmatrix}$ where $d_1=c_3a_2-c_2a_3, d_2=c_1a_3-c_3a_1, d_3= c_2a_1-c_1a_2$. Since $d_1,d_2,d_3\in I$, then $\varepsilon\in E_4(R,I)$. It is easy to check that $\theta(c,b)= \varepsilon^T\theta(a,b)\varepsilon$. Hence $[\theta(c,b)]=[\theta(a,b)]$ in $W_E(R,I)$.

The following theorem is proved in \cite[Theorem 5.2]{SV} in the absolute case. But we need this theorem in the relative case. And the proof is similar to the absolute case.

\begin{theo}
\label{1.7}
 Let $R$ be a commutative ring with $1$ and  $I\subset R$ be an  ideal of $R$. The map $s_E: Um_3(R,I) \longrightarrow W_E(R,I)$ possesses the following properties:

(i) $s_E(b)= s_E(b\alpha)$ for all $b\in Um_3(R,I)$ and $\alpha\in E(R,I)\cap SL_3(R,I)$.

(ii) If $e_1(E(R,I)\cap SL_{2r+1}(R,I))=Um_{2r+1}(R,I)$ for all $r\geq 2$, then $s_E(Um_3(R,I))=W_E(R,I)$.

(iii) If $e_1E_{2r}(R,I)=e_1(E(R,I)\cap SL_{2r}(R,I))$ for all $r\geq 2$, then $s_E(b)=s_E(d)$, where $b,d\in Um_3(R,I)$, only if $b=d\alpha$ for some $\alpha \in E(R,I)\cap SL_3(R,I)$.
\end{theo}

\begin{rem}
 By (i) of Theorem \ref{1.7}, $s_E$ induces a map from $Um_3(R,I)/E_3(R,I)$ to $W_E(R,I)$. We call this map by relative Vaserstein symbol and  denote by $V_{R,I}$.
\end{rem}

\begin{lem}
\label{1.43}
$($\cite[Theorem 5.2]{SV}$)$
Let $A$ be a commutative ring for which $Um_r(A)=e_1E_r(A)$, for $r\geq 5$. Then $V_A$ is surjective. If, moreover, $SL_4(A)\cap E(A)=E_4(A)$, then $V_A$ is bijective. In particular, if $Sd(A) \leq 3$, and $SL_4(A)\cap E_5(A)=E_4(A)$, then $V_A$ is bijective.
\end{lem}

In view of Lemma \ref{1.43}, we can derive the same for relative case from Theorem \ref{1.7}. 

\begin{lem}
\label{1.44}
Let $R$ be a commutative ring and $I\subset R$ be an ideal. If $Um_r(R,I)= e_1E_r(R,I)$ for $r\geq 5$, then $V_{R,I}$ is surjective. If, moreover, $SL_4(R,I)\cap E(R,I)=E_4(R,I)$ then $V_{R,I}$ is bijective. In particular, if $Sd(R)\leq 3$, and $SL_4(R,I)\cap E_5(R,I)=E_4(R,I)$, then $V_{R,I}$ is bijective.
\end{lem}

\begin{cor}
\label{5.10}
Let $R$ be a commutative ring of dimension $2$ and $I\subset R$ be an ideal. Then $Um_3(R,I)/E_3(R,I)$ is bijectively equivalent to $W_E(R,I)$.
\end{cor}

\begin{lem}
\label{4.6}
\cite[Lemma 5.2]{ggr}
Let $R$ be a commutative ring and $I\subset R$ be an ideal of $R$. Let $u,v\in Um_3(R,I)$ be such that $u\alpha=v$ for some $\alpha\in SL_3(R,I)\cap E_4(R,I)$. Then $u$ and $v$ are elementary equivalent relative to $I$.
\end{lem}

\begin{lem}
\label{4.5}
Let $R$ be a commutative ring and $I\subset R$ be an ideal of $R$. Let $\varepsilon\in E_n(R,I)$. Then there exists a matrix $\tilde{\varepsilon}\in E_n(R\oplus I,0\oplus I)$ such that $\varphi(\tilde{\varepsilon})=\varepsilon$, where $\varphi: R\oplus I\rightarrow R$ is defined by $\varphi((r,i))=r+i$.
\end{lem}

\begin{lem}
Let $R$ be a ring and $I\subset R$ be an ideal of  $R$. Assume that both $MSE_3(R)$ and $MSE_3(R\oplus I)$ have Witt group structure via Vaserstein symbol. Let $K$ be the kernel of the map $MSE_3(R\oplus I)\rightarrow MSE_3(R)$ induced from the map $pr_1: R\oplus I \rightarrow R$ defined by $pr_1((r,i))=r$. Then there is a bijection from $MSE_3(R,I)$ to $K$.
\end{lem}
 
 ${\pf}$ Let us define a map $i: MSE_3(R,I) \rightarrow K$ by $i([v])= [\tilde{v}]$, where $v=(1+i_1,i_2,i_3)$ with $i_1,i_2,i_3\in I$ and $\tilde{v}= ((1,i_1),(0,i_2),(0,i_3))$. This map is well defined because if $[v]=[w]$ in $MSE_3(R,I)$, then $v\varepsilon =w$ for some $\varepsilon\in E_3(R,I)$. By Lemma \ref{4.5}, there exists $\tilde{\varepsilon}\in E_3(R\oplus I, 0\oplus I)$ such that $\varphi(\tilde{\varepsilon})= \varepsilon$. It is easy to check that $\tilde{v}\tilde{\varepsilon}=\tilde{w}$.
 
 \begin{lem}
 Let $R$ be a commutative ring and $I\subset R$ be an ideal of $R$. Then the relative Witt groups  $W_E(R,I)$ and $W_E(R\oplus I,0\oplus I)$ are isomorphic.
 \end{lem}
  
  ${\pf}$
  Let us define a map $\varphi: W_E(R,I)\rightarrow W_E(R\oplus I, 0\oplus I)$ by $\varphi([\alpha])=[\tilde{\alpha}]$. Clearly this map is well-defined. Let $\varphi([\alpha])=[\chi_1]$. Then there exists $t\in \mathbb{N}$ and  $\varepsilon\in E_{2(n+t)}(R\oplus I, 0\oplus I)$ such that $\varepsilon^T(\tilde{\alpha}\perp \chi_t)\varepsilon= \chi_{n+t}$. Projecting $R\oplus I$ onto $R$ we have $\varepsilon_1^T(\alpha\perp \chi_t)\varepsilon_1= \chi_{n+t}$, where $\varepsilon\in E_{2(n+t)}(R,I)$. Hence $\varphi$ is injective. Clearly $\varphi$ is surjective. It is easy to check that $\varphi$ is a group homomorphism. Hence $\varphi$ is an isomorphism.
  $~~~~~~~~~~~~~~\qedwhite$

   Following L.N. Vaserstein in \cite{SV}, one can show the following: 

\begin{lem} $($Relative version of L.N. Vaserstein's lemma$)$\label{re2}

Let $\varphi \in {\rm SL}_{2n}(R)$ be an alternating matrix. Then 
\begin{eqnarray*}
e_1{\rm E}_{2n}(R, I) &=& e_1{\rm Sp}_{\varphi}(R, I) \cap {\rm E}_{2n}(R, I).
\end{eqnarray*} 

Here we have the isotropy group  
\begin{eqnarray*}
{Sp}_{\varphi}(R, I) & = &
\{\varepsilon \in {SL}_{2n}(R, I)|~\varepsilon^t \varphi \varepsilon = 
\varphi\}.
\end{eqnarray*}
\end{lem}

${\pf}$
Let $\varepsilon \in {\rm E}_{2n}(A, I)$, and $\varepsilon(X) \in 
{\rm E}_{2n}(A[X], I[X])$ be a homotopy of $\varepsilon$. Let $v(X)= e_1
\varepsilon(X)$. Let $\mathfrak{m}$ be a maximal ideal of $A$. By 
(\cite{SV}, Lemma 5.5)  $v(X)_\mathfrak{m} = e_1 \varepsilon_\mathfrak{m}(X)$, 
for some 
$\varepsilon_\mathfrak{m}(X)$, with $\varepsilon_\mathfrak{m}(0) = I_{2n}$, 
and with 
\begin{eqnarray*}
\varepsilon_\mathfrak{m}(X) &\in & 
\{{Sp}_\varphi)_{2n}(A_\mathfrak{m}[X], I_\mathfrak{m}[X]) \cap 
 {E}_{2n}(A_\mathfrak{m}[X], I_\mathfrak{m}[X])\}
\end{eqnarray*}
By imitating the proof of the action version of the Local Global Theorem \cite[Theorem 4.7]{acr}, one can show that  there is a 
global $\varepsilon_\varphi(X)$ with $v(X) = e_1 \varepsilon_\varphi(X)$, 
with $\varepsilon_\varphi(X) \in 
{Sp}_\varphi(A[X], I[X]) \cap 
 {E}_{2n}(A[X], I[X])\}$. Now put $X = 1$. 

\vskip0.15in

In fact, Chattopadhyay--Rao have shown that

\vskip0.15in

\begin{theo} $($\cite{cr}$)$ \label{cr}

Let $R$ be a commutative ring with $1$. Let $I$ be an ideal of $R$. Let 
$n \geq 2$, and let  $v \in Um_{2n}(R, I)$. Then 
$vE_{2n}(R, I) = vESp_{2n}(R, I)$. 
\end{theo}

\vskip0.15in

\vskip0.15in

We state a lemma which can be found in (\cite{MG}, Lemma 4.3) where it is 
attributed to A. Suslin: 

\vskip0.15in

\begin{lem} \label{k2}
Let $A$ be any commutative ring with $1$. Let $I$ be an ideal 
of $A$ such that {\rm K}$_{2, n}(A) \longrightarrow {\rm K}_{2, n}(A/I)$, 
for some $n \geq 3$, is surjective. Then ${\rm E}_n(A, I) = {\rm E}_n(A) 
\cap {\rm GL}_n(A, I)$. \hfill$\Box$
\end{lem}

\vskip0.15in

\begin{theo}\label{relvas}

Let $R$ be a commutative ring of dimension $2$. Then the Vaserstein 
symbol 
\begin{eqnarray*}
\frac{{\rm Um}_3(R[T], (T^2 - T))}{{\rm E}_3(R[T], (T^2- T))} & 
\longrightarrow & {\rm W}_{\rm E(R[T], (T^2 - T))}
(R[T], (T^2 - T))
\end{eqnarray*}
is surjective. If $R$ is an affine algebra over a perfect {\rm C}$_1$-field  
then $V$ is also injective if $6R = R$. 
\end{theo}

${\pf}$ Let $A = R[T]$, $t = T^2 - T$, $x = X^2 - X$ . The surjectivity 
will follow by $(ii)$ of Theorem \ref{1.7}.

For the injective, the argument of L.N. Vaserstein in \cite{SV} says, in 
view of this, that it suffices to 
show that if $\sigma \in {\rm Sl}_4(A, (t)) \cap {\rm E}(A, (t))$ 
then $\sigma \in {\rm E}_4(A, (t))$. By stability estimates, 
$\sigma \in {\rm E}_5(A, (t))$. 

Choose a homotopy $\rho(X) \in {\rm E}_5(A[X], (tX)])$ of $(1 \perp \sigma)$. 
Clearly,  $e_1 \rho(X) \in {Um}_5(A[X], (tx))$. By (\cite{raostably}, 
Proposition 5.4)  there is a $\varepsilon(X) \in {\rm E}_5(A[X], 
(tx))$ such that $e_1 \rho(X) \varepsilon(X)$ is a factorial 
vector; whence it can be completed to a $\beta(X)$ which is stably 
elementary and congruent to identity modulo $(tx)$ (see 
(\cite{raostably}, Remark 5.5)). If 
\begin{eqnarray*}
\rho(X) \varepsilon(X) \beta(X)^{-1} & = & \begin{pmatrix} 1 & 0 \\
* & \rho^*(X)\end{pmatrix}, 
\end{eqnarray*}
then $\rho^*(X)$ is a (stably elementary) homotopy of $\sigma$. By 
(\cite{S}, Theorem 6.3), $\rho^*(X)$ is an elementary matrix. 

For any maximal ideal $\mathfrak{m}$ of $A$, we have $\rho^*(X)_{\mathfrak{m}}\in SL_4(A_{\mathfrak{m}}[X],(t)\cap E_4(A_{\mathfrak{m}}[X])$. Hence by Lemma \ref{k2}, we have $\rho^*(X)_\mathfrak{m}\in E_4(A_{\mathfrak{m}}[X],(t))$.  Hence, by the Local-Global theorem for an extended ideal 
\cite[Theorem 4.5]{acr}, we have  $\rho^*(X) \in {\rm E}_4(A[X], (t))$. 

Substituting $X = 1$, gives $\sigma \in {\rm E}_4(A, (t))$.
$~~~~~~~~~~~~~~~~~~\qedwhite$

An alternative approach to Theorem \ref{relvas} in case $R$ is a non-singular affine algebra of dimension $2$ over a perfect $C_1$-field is as follows:

  \begin{theo}
  \label{dim3principle}
  \cite[Theorem 4.7]{anjan}
  Let $R$ be an affine non-singular algebra of dimension $d$ over a perfect $C_1$-field $k$, and $I\subset R$ be a principal ideal of $R$. Let $\alpha \in SL_{d+1}(R,I)\cap E(R,I)$, then $\alpha\in E_{d+1}(R,I)$.
  \end{theo}

\begin{cor}
Let $R$ be a non-singular affine algebra  of dimension $2$ over a perfect $C_1$-field $k$. Then the Vaserstein symbol $$V_{R[X],(X^2-X)}:MSE_3(R[X],(X^2-X))\rightarrow W_E(R[X],(X^2-X))$$ is bijective.
\end{cor}

${\pf}$
The surjectivity is clear from $(ii)$ of Theorem \ref{1.7}. In view of $(iii)$ of Theorem \ref{1.7}, $V$ is injective if we can show that $SL_4(R[X],(X^2-X))\cap E(R[X],(X^2-X))= E_4(R[X],(X^2-X))$. Which follows from Theorem \ref{dim3principle}.
$~~~~~~~~~~~\qedwhite$

\begin{theo}
\label{1.22.0}
$($\cite[Theorem 4.8]{anjan}$)$
 Let $A$ be a non singular affine algebra of dimension $d$ over an algebraically closed field $k$, $d\geq 3$, $d! \in k^*$ and $I=(a)$ a principal ideal of $A$. Then
$SL_d(A, I)\cap  E(A, I) = E_d(A, I)$.
\end{theo}

\begin{cor}
Let $R$ be a non-singular affine algebra of dimension $3$ over an algebraically closed field $k$, $4!\in k^*$. Then the map $V_{R[X],(X^2-X)}: MSE_3(R[X],(X^2-X))\rightarrow W_E(R[X],(X^2-X))$ is injective. 
\end{cor}

${\pf}$ Start with the map $s_E: Um_3(R[X],(X^2-X))\rightarrow W_E(R[X],(X^2-X))$. Let $s_E(b)=s_E(d)$, for some $b,d\in Um_3(R[X],(X^2-X))$. In view of Theorem \ref{1.7}, the map $s_E$ is injective if $b=d\alpha$ for some $\alpha\in SL_3(R[X],(X^2-X))\cap E(R[X],(X^2-X))$. By Theorem \ref{1.22.0}, we have $\alpha\in SL_3(R[X],(X^2-X))\cap E_4(R[X],(X^2-X))$. But by Lemma \ref{4.6}, we have $b,d$ are elementary equivalent relative to $(X^2-X)$. Hence the map $s_E$ induces an injective map $V_{R[X],(X^2-X)}: MSE_3(R[X],(X^2-X))\rightarrow W_E(R[X],(X^2-X))$.
$~~~~~~~~~~~~~~~\qedwhite$

\begin{theo}
\cite[Theorem 4.1]{anjan}
\label{relsuslin}
Let $R$ be an affine algebra of dimension $d$, $d\geq 2$, over a perfect $C_1$-field $k$. Then $Um_{d+1}(R,I)=e_1SL_{d+1}(R,I)$, for any ideal $I$ of $R$.
\end{theo}

\begin{theo}
\label{localcompinte}
Let $R$ be an affine non-singular algebra of dimension $3$ over  a  perfect $C_1$-field $k$ and $I\subset R$ be a local complete intersection ideal of $R$. Then the map $V_{R,I}: MSE_3(R,I)\rightarrow W_E(R,I)$ is injective.
\end{theo}

${\pf}$ In view of Lemma \ref{1.44}, it is enough to show that $SL_4(R,I)\cap E(R,I)= E_4(R,I)$. Let $\alpha\in SL_4(R,I)\cap E(R,I)$. By classical stability estimates, we have $1\perp \alpha\in E_5(R,I)$. Hence there exists $\alpha(X)\in E_5(R[X],I[X])$ such that $\alpha(0)= I_5$ and $\alpha(1)=1\perp \alpha$. Let $v(X)= e_1\alpha(X)$. We have $v(X)\in Um_5(R[X], (X^2-X)I[X])$. Since $R[X]$ is affine algebra of dimension $4$ over perfect $C_1$-field $k$, we have $v(X)\in e_1SL_5(R[X],(X^2-X)I[X])$ by Theorem \ref{relsuslin}. By \cite[Proposition 5.4]{raostably}, we have there exists $\varepsilon(X)\in E_5(R[X],(X^2-X)I[X])$ and $\beta(X)\in SL_5(R[X],(X^2-X)I[X])\cap E(R[X],(X^2-X)I[X])$, such that $v(X)\varepsilon(X)\beta(X)=e_1$. Hence by arguing as in \cite[Theorem 3.4]{rv}, we have a relative stably elementary homotopy for $\alpha$. Hence there exists $\alpha(X)\in SL_4(R[X],I[X])\cap E_5(R[X],I[X])$. By Local-Global principle for an extended ideal, we may assume that $R$ is a local ring with the maximal ideal $\mathfrak{m}$. Let $S= R\setminus\{0\}$. Then we have $\alpha(X)_{S}\in SL_4(R_S[X],I_S[X])\cap E_5(R_S[X],I_S[X])= E_4(R_S[X],I_S[X])$. Hence there exists $h\in R$, such that $\alpha(X)_h\in E_4(R_h[X],I_h[X])$. We may assume that $h\in \mathfrak{m}$. Further we may assume that $h\in \mathfrak{m}^2$.  By  Theorem \ref{1.36}, we have, there exists $L$, a subring of $R$ such that 

(i) $L\hookrightarrow R$ is analytic isomorphism along $h^\prime$ for some $h^\prime\in hR\cap L $.

(ii) $L=K[X_1,X_2,X_3]_{(\varphi(X_1),X_2,X_3)}$, where $\varphi(X_1)\in K[X_1]$ is a monic irreducible polynomial.

By Lemma \ref{1.37}, we have $h$ and $h^\prime$ differ by a unit in $R$, hence therefore we may assume that $h=h^\prime$. Therefore we have , $L[X]\hookrightarrow R[X]$ is analytic isomorphism along $h$. 

Further we may assume that the ideal $I$ is an ideal of $L$ and $I$ is generated by a subset of the set $\{X_2,X_3\}$. This is possible since $I$ is a local complete intersection. 
Also by Lemma \ref{1.38}, we have $L[X]\oplus I[X]\hookrightarrow R[X]\oplus I[X]$ is analytic isomorphism along $(h,0)$.

Let $\tilde{\alpha(X)}$ be a lift of $\alpha(X)$ in $SL_4(R[X]\oplus I[X])$. Also we have $\tilde{\alpha(X)}_{(h,0)}\in E_4(R_h[X]\oplus I_h[X])$, i.e., $\tilde{\alpha(X)}_{(h,0)}\in E_4((R[X]\oplus I[X])_{(h,0)})$. Then by Lemma \ref{1.37}, we have there exists $\tilde{\gamma(X)}\in E_4(R[X]\oplus I[X])$ and $\tilde{\beta(X)}\in SL_4(L[X]\oplus I[X])$ such that $\tilde{\alpha(X)}=\tilde{\gamma(X)}\tilde{\beta(X)}$.

Now going modulo $0\oplus I[X]$, we have $I_4=\bar{\tilde{\gamma(X)}}\bar{\tilde{\beta(X)}}$, where $\bar{\tilde{\gamma(X)}}\in E_4(R[X])\subseteq E_4(R[X]\oplus I[X])$ and $\bar{\tilde{\beta}}\in SL_4(R[X])\subseteq SL_4(R[X]\oplus I[X])$. Hence replacing $\tilde{\gamma(X)}$ by $\tilde{\gamma(X)}\bar{\tilde{\gamma(X)}}^{-1}$ and $\tilde{\beta(X)}$ by $\bar{\tilde{\beta(X)}}^{-1}\tilde{\beta(X)}$, we may assume that $\tilde{\gamma(X)}\in SL_4(R[X]\oplus I[X],0\oplus I[X])\cap E_4(R[X]\oplus I[X])=E_4(R[X]\oplus I[X])$ and $\tilde{\beta(X)}\in SL_4(L[X]\oplus I[X],0\oplus I[X])$.

Hence projecting $R[X]\oplus I[X]$ onto $R[X]$ and $L[X]\oplus I[X]$ onto $L[X]$ we have $\alpha(X)=\gamma(X)\beta(X)$, where $\gamma(X)\in E_4(R[X],I[X])$ and $\beta(X) \in SL_4(L[X],I[X])$. Also  we have $\beta(X)_h\in E_4(L_h[X], I_h[X])$.

We may assume that $h$ belongs to the maximal ideal of $L$. Also by multiplying unit of $L$ we may assume that $h\in (\phi(X_1),X_2,X_3)$. Now by Lemma \ref{3.5}, we have, there is a transformation of $K[X_1,X_2,X_3]$, namely $X_1\mapsto X_1, X_i\mapsto X_i+\phi(X_1)^{r_i}, i\geq 2$, for some $r_i$ such that $h$ becomes a monic polynomial in $X_1$ over $K[X_2,X_3]$.

Let $L^{\prime}=K[X_2,X_3]_{(X_2,X_3)}[X_1]$. Since the above transformation takes the maximal ideal $(\phi(X_1),X_2,X_3)$ of $K[X_1,X_2,X_3]$ to itself, then, the polynomial $h\in L^{\prime}$ is a Weierstrass polynomial. Hence  by Lemma \ref{0.6}, we have $L^{\prime}\hookrightarrow L$ is analytic isomorphism along $h$.

Hence we have $L^{\prime}[X]\hookrightarrow L[X]$ is analytic isomorphism along $h$. Since $\beta(X)_h\in E_4(L_h[X],I_h[X])$, then again by  Lemma \ref{1.37}, there exists $\delta(X) \in E_4(L[X], I[X])$  and $\theta \in SL_4(L^{\prime}[X], I[X])$ such that $\beta(X)=\delta(X)\theta(X)$ and $\theta(X)_h\in E_4(L^{\prime}_h[X],I_h[X])$. Since $I[X]$ is an extended ideal of $L^{\prime}[X]= (k[X_2,X_3]_{(X_2,X_3)}[X])[X_1]$, and $h$ is monic polynomial in $X_1$ of $L^{\prime}[X]$, then by Lemma \ref{3.4}, we have, $\theta\in E_4(L^{\prime}[X],I[X])$.

Hence $\beta(X)\in E_4(L[X],I[X])\subseteq E_4(R[X], I[X])$. Finally we have $\alpha(X) \in E_4(R[X], I[X])$.

Hence by Local-Global principle for an extended ideal, we have $\alpha(X)\in E_4(R[X],I[X])$, where $R$ is affine non-singular algebra of dimension $3$ over perfect $C_1$-field $k$. Hence evaluating at $X=0$, we have $\alpha \in E_4(R,I)$. 

Hence the  map $V_{R,I}$ is injective.

$~~~~~~~~~~~~~~~~~~~~~~~~~~~~~~~~~~~~~~~~~~~\qedwhite$

\begin{theo}
\cite[Theorem 4.4]{anjan}
\label{acfsuslin}
Let $R$ be a non-singular affine algebra of dimension $d\geq 4$ over an algebraically closed field $k$, $1/(d-1)! \in k$ and $I\subset R$ be an ideal of $R$. Let $v\in Um_d(R,I)$. Then $v$ can be transformed to a factorial row $\psi_{(d-1)!}(w)=(w_0,w_1, w_2^2\dots,w_{d-1}^{d-1})$ for some $w=(w_0,w_1,\dots,w_{d-1})\in Um_d(R, I)$ by elementary operations relative to $I$. In particular $Um_d(R,I)=e_1SL_d(R, I)$.
\end{theo}

The proof of the following theorem is same as the proof of Theorem \ref{localcompinte}, only we have to use Theorem \ref{acfsuslin}, instead of Theorem \ref{relsuslin}.

\begin{theo}
\label{localcomintedim4}
Let $R$ be a non-singular affine algebra of dimension $4$ over an algebraically closed field $k$, and let $I\subset R$ be a local complete intersection ideal of $R$. Then the Vaserstein symbol $$V_{R,I}: MSE_3(R,I)\rightarrow W_E(R,I)$$ is injective. 
\end{theo}

The proof of the following theorem is same as the proof of Theorem \ref{localcompinte}.

\begin{theo}
\label{extendedvas}
Let $R$ be a non-singular affine algebra of dimension $2$ over a perfect $C_1$-field $k$, and let $I\subset R$ be a local complete intersection ideal of $R$. Then the Vaserstein symbol $$V_{R[X],I[X]}:MSE_3(R[X],I[X])\rightarrow W_E(R[X],I[X]))$$ is injective.
\end{theo}

\section{\bf{Kernel of the Vaserstein symbol}}

   \begin{theo}
   \cite[Theorem 3.6]{vdk1}
   Let $R$ be a noetherian ring of dimension $d$, and $I\subset R$ be an ideal of $R$. Then $MSE_{d+1}(R,I)$ is  an abelian group with the following operation:
   
If $[v], [w]\in MSE_{d+1} (R,I)$, choose representatives $(a_0,a_1,\dots,a_d)\in [v], (b_0,b_1,\dots,b_d)\in [w]$ with $a_i=b_i$ for $i\geq 1$, and choose $p_0\in R$ such that $a_0p_0\equiv 1 \mod (a_1R+a_2R+\dots +a_dR)$. Then 
$$[w]. [v]= [(a_0(b_0+p_0)-1,(b_0+p_0)a_1,a_2,\dots,a_d)].$$
   \end{theo}

\begin{defi}
Let $R$ be a ring and $I\subset R$ be an ideal of $R$. Then excision ring $\mathbb{Z}\oplus I$ is defined by the set $\{(n,i): n\in \mathbb{Z}, i\in I\}$ with the following operations:

\begin{enumerate}
\item $(m,i)+(n,j)=(m+n,i+j)$
\item $(m,i)(n,j)=(mn,mj+ni+ij)$
\end{enumerate}
\end{defi}

\begin{theo}[Excision theorem]
\label{excisiontheo}
\cite[Theorem 3.21]{vdk1}
Let $R$ be a ring and $I\subset R$ be an ideal od $R$. Then for $n\geq 3$, the natural maps $MSE_n(\mathbb{Z}\oplus I, 0\oplus I)\overset{f}\rightarrow MSE_n(R,I)$ defined by $f(mse(v))= mse(f(v))$ and $MSE_n(\mathbb{Z}\oplus I, 0\oplus I)\overset{i}\rightarrow MSE_n(\mathbb{Z}\oplus I)$ defined by $i(mse(v))=mse(v)$ are bijective. 
\end{theo}

\begin{theo}[van der Kallen]
\label{groupstr}
\cite[Theorem 4.1]{vdk2}
Let $R$ be a ring of stable dimension $d$, $d\leq 2n-4$, $n\geq 3$. Then the universal weak Mennicke symbol $wms: MSE_n(R)\rightarrow WMS_n(R)$ is bijective and hence $MSE_n(R)$ has an abelian group structure.
\end{theo}

\begin{rem}
Let $R$ be a ring of stable dimension $d$, $d\geq 2$ and $I$ be an ideal of $R$. Then by \cite[3.19]{vdk1}, we have the maximal spectrum of $\mathbb{Z}\oplus I$ is the union of finitely many subspaces of dimension $\leq d$. Therefore by Theorem \ref{groupstr}, we have $MSE_n(\mathbb{Z}\oplus I)$ has a group structure for $n\geq \max\{3,d/2+2\}$. Hence we can say the group structure of $MSE_n(R,I)$ for $n\geq \max\{3,d/2+2\}$ via the excision theorem \ref{excisiontheo}.
\end{rem} 

\begin{defi}
Let $R$ be a ring with $I\subset R$ be an ideal of $R$. We call the group structure of $MSE_n(R)$ given by van der Kallen (Theorem \ref{groupstr}) is nice if it satisfies the following  'coordinate-wise multiplication' formula:
$$[(b,a_2,\dots, a_n)].[(a,a_2,\dots, a_n)]=[(ab,a_2,\dots,a_n)].$$
 
 Similarly we call the group structure of $MSE_n(R,I)$ given by van der Kallen is nice if it satisfies the following 'coordinate-wise multiplication' formula:
 $$[(b,a_2,\dots, a_n)].[(a,a_2,\dots, a_n)]=[(ab,a_2,\dots,a_n)].$$
\end{defi}

 \begin{theo}
 \label{4.10.1}
 \cite[Theorem 2.2]{vdk2}
 Let $R$ be a commutative ring with $sdim(R)\leq 2n-3$, $n\geq 3$ and $I$ be an ideal of $R$. Let $i,j$ be non-negative integers. Then for every $g\in GL_{n+i}(R)\cap E_{n+i+j+1}(R,I)$, there exist matrices $u,v,w,M$ with entries in $I$ and $q$ with entries in $R$ such that \begin{center}
$ \begin{pmatrix}
 I_{i+1}+uq &v\\
 wq &I_{n-1}+M
 \end{pmatrix} \in gE_{n+i}(R,I)$,
 
 $\begin{pmatrix}
 I_{j+1}+qu &qv\\
 w&I_{n-1}+M
 \end{pmatrix}\in E_{n+j}(R,I)$.
 \end{center}
  \end{theo}

\begin{cor}
\label{niceimplytrivial}
Let $R$ be a commutative ring of dimension $d$ with orbit space $MSE_n(R,I)$, $n\geq \max\{3,d/2+2\}$ has a nice group structure. Let $\sigma\in SL_n(R,I)\cap E_{n+1}(R,I)$, then $[e_1\sigma]=1$ in $MSE_n(R,I)$.
\end{cor}

${\pf}$ Putting $i,j=0$ in the Theorem \ref{4.10.1}, we have,
\begin{center}
$ \begin{pmatrix}
 1+uq &v\\
 wq &I_d+M
 \end{pmatrix} \in \sigma E_n(R,I)$,
 
 $\begin{pmatrix}
 1+qu &qv\\
 w&I_d+M
 \end{pmatrix}\in E_n(R,I)$.
 \end{center}
 Hence in the orbit space $MSE_n(R,I)$ we have 
 \begin{align*}
 [e_1\sigma]\\
 &=[(1+uq,v)]\\
 &=[(1+uq,v)][1+uq,q];~~~~\text{since $[(1+uq,q)]$ is the  identity}\\
 &=[(1+uq,vq)];~~~~\text{since $MSE_n$ has nice group structure}\\
 &=[e_1].
 \end{align*}
$~~~~~~~~~~~~~~~~~~~~~~~~~~~~~~~\qedwhite$

\begin{theo}
\label{kernel}
Let $R$ be a ring of dimension $3$ and $I\subset R$ be an ideal of $R$. If the orbit space $MSE_4(R,I)$ has a nice group structure, then,  the Vaserstein map $V_{R,I}$ induces a bijection  of the map $$\phi: Um_3(R,I)/\{\sigma\in SL_3(R,I)\cap E_5(R,I)\}\equiv W_E(R,I).$$
\end{theo}

${\pf}$ Since dimension of $R$ is $3$, then $\phi$ is surjective by Theorem \ref{1.7}. Let $v\in ker(\phi)$. Then there exists $\tau\in SL_4(R,I)\cap E_5(R,I)$ such that $\tau^TV(v,w)\tau=\chi_2$.  Since $MSE_4(R,I)$ has a nice group structure, then, by Corollary \ref{niceimplytrivial}, we have $e_1\tau$ is elementarily completable. Hence we have $\tau\varepsilon=1\perp \rho$, with $\rho\in SL_3(R,I)\cap E_5(R,I)$ and $\varepsilon\in E_4(R,I)$. Now from \cite[Chapter 5]{SV}, we have $v\rho \varepsilon_0=e_1$ for some $\varepsilon_0\in E_3(R,I)$ as desired.
$~~~~~~~~~~~~~~\qedwhite$

\begin{theo}
\cite[Theorem 4.2]{ggr}
Let $R$ be an affine algebra of dimension $d$ over a perfect $C_1$-field $k$, with char $k \neq 2$. Let $I\subset R$
 be an ideal of $R$. Then the group structure of $MSE_{d+1}(R,I)$ is nice.
 \end{theo}
 
 \begin{cor}
 Let $R$ be an affine algebra of dimension $3$ over a perfect $C_1$- field $k$ with char $k \neq 2$, then the Vaserstein map $V_{R,I}$ induces a bijection  $$\phi: Um_3(R,I)/\{\sigma\in SL_3(R,I)\cap E_5(R,I)\} \equiv W_E(R,I).$$
 \end{cor}

\section{\bf{The Vaserstein symbol for Rees algebras and extended Rees algebras}}

\begin{defi}
${\bf{Rees ~ algebra~and~extended~Rees~algebras:}}$
Let $R$ be a commutative ring of dimension $d$ and $I$ an ideal of $R$. Then the algebra
\[
R[It] :=\{\sum_{i=0}^{n}a_it^i :n\in \mathbb{N}, a_i\in I^i\}= \oplus_{n\geq 0}I^nt^n
\]
is called the Rees algebra of $R$ with respect to $I$.
The extended Rees algebra of $R$ with respect to $I$, denoted by  $R[It,t^{-1}]$,is defined by
\[
R[It,t^{-1}] := \{\sum_{i=-n}^{n}a_it^i : n\in \mathbb{N}, a_i\in I^i\} = \oplus_{n\in \mathbb{Z}} I^nt^n
\]
where $I^n=R$ for $n\leq 0$.
\end{defi}

Clearly the Rees algebra $R[It]$ is a graded ring. The following graded version of Quillen's Local-Global Principle is well-known: 

\begin{theo}
\label{2.0}
$($\cite[Theorem 4.3.11]{ir}$)$
Let $S=S_0\oplus S_1\oplus S_2\oplus ...$ be a commutative graded ring and let $M$ be a finitely presented S-module. Assume that for every maximal ideal $\mathfrak{m}$ of $S_0$,  $M_{\mathfrak{m}}$ is extended from $(S_0)_{\mathfrak{m}}$. Then $M$ is extended from $S_0$.
\end{theo}

\begin{lem}
\label{2.1}
$($\cite [Lemma 3.1] {rs}$)$
 Let $R$ be a commutative ring and $I,J$ ideals of $R$. Then the natural map $\phi : R[It]/JR[It] \rightarrow \bar R[\bar It]$, where $\bar{R} = R/J$, $\bar I = (I+J)/J$, defined by $\phi( r+ a_1t+a_2t^2 +\dots+a_nt^n +JR[It]) =  \bar r+\bar a_1t+\dots+\bar  a_nt^n$, is an isomorphism.
\end{lem}

\begin{theo}
\label{2.2}
$($\cite[Theorem 1.3]{Vasco}$)$
Let $R$ be a Noetherian ring of dimension $d$ and $I$ be an ideal of $R$. Then dimension of $R[It]\leq d+1$. Moreover if $I$ is not contained in any minimal primes of $R$, then $dim(R[It])=d+1$.
\end{theo}

\begin{theo}
\label{2.3}
$($\cite[Theorem 2.1]{rs}$)$
Let $R$ be a Noetherian regular ring. Then $R[It]$ is regular if and only if $I=(0)$ or $(1)$ or generated by single element.

\end{theo}

\begin{theo}
\label{2.4}
$($\cite[Theorem 4.2]{rs}$)$
Let $R$ be a ring of dimension $d$ and $I\subset R$ an ideal of $R$. Then for $n\geq max\{3, d+2\}$, the natural map $\phi : GL_n(R[It])/E_n(R[It])\rightarrow K_1(R[It])$ is an isomorphism.
\end{theo}

\begin{theo}
\cite[Theorem 5.1]{overring}
\label{overring}
Let $R$ be a commutative ring of dimension $d$ and $A$ be a ring lying between $R[X]$ and $S^{-1}R[X]$, where $S$ is a multiplically closed set of non-zero-divisors in $R[X]$, then for $n\geq max\{3, d+2\}$, the natural map $\phi: GL_n(A)/E_n(A)\rightarrow K_1(A)$ is an isomorphism.

\end{theo}

\begin{cor}
\label{ex-stabi}
Let $R$ be a ring of dimension $d$ and $I\subset R$ be an ideal of $R$. Then for $n\geq \max\{3, d+2\}$, the natural map $\phi: GL_n(R[It,t^{-1}])/E_n(R[It,t^{-1}]) \rightarrow K_1(R[It,t^{-1}])$ is an isomorphism.
\end{cor}

${\pf}$ We have $R[t^{-1}]\subset R[It,t^{-1}]\subset R(t^{-1})$, where $R(t^{-1})$ is the total quotient ring of $R[t^{-1}]$. Now apply Theorem \ref{overring}.

Here we are giving an example of a  $3$-dimensional algebra  for which the Vaserstein symbol is  injective though the algebra is singular.

\begin{theo}
\label{2.6}
Let $R$ be an affine algebra of dimension $2$  and $I\subset R$ be an ideal. Then the Vaserstein symbol $V_{R[It]}: MSE_3(R[It])\rightarrow W_E(R[It])$ is bijective.
\end{theo}

${\pf}$ By Theorem \ref{2.2} we have dimension of $R[It]\leq 3$. Hence we have $Sd(R[It])\leq 3$.
By Theorem \ref{1.43}, it is enough to show that $SL_4(R[It])\cap E_5(R[It])= E_4(R[It])$.
 By Theorem \ref{2.4}, we have for $n\geq 4$, the natural map $\phi : GL_n(R[It])/E_n(R[It])\rightarrow K_1(R[It])$ is an isomorphism. So in particular for $n=4$, we have  $SL_4(R[It])\cap E_5(R[It])= E_4(R[It])$.
 $~~~~~~~~~~~~~~~~~~~~~~~~~~~~~~~~~~~~~~~~~~~~~~~~~~~~~~~~
 \qedwhite$

\begin{notn}
Let $R$ be a commutative ring and $I\subset R$ be an ideal. Let $S\subset R $ be a subring of $R$. Then we denote the algebra $S[It]$ by 
\[
S[It] := \{\sum_{i=0}^{n}a_it^i :n\in \mathbb{N},a_i\in I^i \} = \oplus_{n \geq 0} I^nt^n
\]
where $I^0=S$.
\end{notn}

$S[It]$ is an algebra since the subring $S$ induces an $S$ module structure on $R$. Hence $I$ is an $S$ module.

\begin{lem}
\label{2.8}
Let $R$ be a commutative ring and $I\subset R$ be an ideal. Let $S \subset R$ be a subring of $R$. If there exists some $h\in S$ for which  $S\hookrightarrow R$ is an analytic isomorphism along $h$, then $S[It]\hookrightarrow R[It]$ is an analytic isomorphism along $h$.
\end{lem}

${\pf}$
 Clearly $h$ is a non-zero divisor of $S[It]$ as well as of  $R[It]$. Since $S \hookrightarrow R$ is an analytic isomorphism along $h$, then the natural map $i: S\rightarrow R$ induces an isomorphism between $S/hS$ and $R/hR$. Therefore the induced map $\bar{i}: S/hS \rightarrow R/hR$ induces an isomorphism between $\bar{S}[\bar{I}t]$ and $\bar{R}[It]$, where $\bar{S}=S/hS$, $\bar{R}=R/hR$ and $\bar{I}=(I+hR)/hR$. Consider the following diagram

\[
\begin{tikzcd}
S[It]/hS[It] \arrow{r}{\bar{i^*}} \arrow{d}{\psi} &R[It]/hR[It]\arrow{d}{\phi}\\
\bar{S}[\bar{I}t] \arrow{r}{\bar{i}^*} &\bar{R}[\bar{I}t]
\end{tikzcd}
\]
The diagram is commutative. 
By Lemma \ref{2.1} we have the map $\phi$ is an isomorphism. It can be shown that the map $\psi$ is also an isomorphism. Therefore by commutativity of the above diagram we have $\bar{i^*}$ is an isomorphism.
$~~~~~~~~~~~~~~~~~~~~~~~~~~~~~~~~~~~~~~~~~~~~~~~~\qedwhite$

The same result holds for extended Rees algebra and the proof is same as Lemma \ref{2.8} has.

\begin{lem}
\label{2.9}
Let $R$ be a commutative ring and $I\subset R$ be an ideal. Let $S\subset R$ be a subring of $R$. If there exists some $h\in R$ for which $S\hookrightarrow R$ is an analytic isomorphism along $h$, then $S[It,t^{-1}]\hookrightarrow R[It,t^{-1}]$ is an analytic isomorphism along $h$.
\end{lem}

\begin{theo}
\label{4.10}
Let $R$ be a non-singular affine algebra of dimension $3$ over a field $k$ and $I\subset R$ be an ideal. Then the Vaserstein symbol $V_{R[It]}: MSE_3(R[It])\rightarrow W_E(R[It])$ is injective.
\end{theo}

${\pf}$
By Lemma \ref{1.43}, it is enough to check that $SL_4(R[It])\cap E_5(R[It])=E_4(R[It])$. Let $\alpha\in SL_4(R[It])\cap E_5(R[It])$. By Local Global principle (Theorem \ref{2.0}) we may assume that, $R$ is a regular local ring. First we  assume that $R$ is regular local algebra with a separating ground field $K$. Let $S=R\setminus \{0\}$. By stability we  have $\alpha_S\in E_4(S^{-1}R[It])$. Hence there exists $g\in R$ such that $\alpha_h\in E_4(R_g[It])$. We may assume that $g\in \mathfrak{m}$, the maximal ideal of $R$. Further we may assume that $g\in \mathfrak{m}^2$. By Theorem \ref{1.36}, there exists a subring $L$ of $R$ and an element $h\in L\cap gR$ such that $L= K[X_1,X_2,\dots,X_d]_{(\varphi(X_1),X_2,\dots,X_d)}$, where $\varphi(X_1)$ is an irreducible monic polynomial, and $L\hookrightarrow R$ is an analytic isomorphism along $h$. By Lemma \ref{2.8}, $L[It]\hookrightarrow R[It]$ is an analytic isomorphism along $h$. Thus we have a patching diagram
\[
\begin{tikzcd}
L[It]\arrow{r}{}\arrow{d}{} &R[It]\arrow{d}{}\\
L_h[It] \arrow{r}{} &R_h[It]
\end{tikzcd}
\]
By Theorem \ref{1.36}, we have $g$ and $h$ are differed by a unit in $R$. Hence we have $\alpha_h\in E_4(R_h[It])$. Hence by Lemma \ref{1.37} we have there exists $\beta\in SL_4(L[It])$ with $\beta_h\in E_4(L_h[It])$ and $\gamma \in E_4(R[It])$ such that $\alpha= \gamma\beta$.  Therefore it is enough to show that $\beta \in E_4(L[It])$.
 
 We may assume that $h$ belongs to the maximal ideal of $L$. Also by multiplying unit of $L$ we may assume that $h\in (\phi(X_1),X_2,...,X_d)$. Now by Lemma \ref{3.5}, we have, there is a transformation of $K[X_1,X_2,...,X_d]$, namely $X_1\mapsto X_1, X_i\mapsto X_i+\phi(X_1)^{r_i}, i\geq 2$, for some $r_i$ such that $h$ becomes a monic polynomial in $X_1$ over $K[X_2,X_3,...,X_d]$.

Let $L^{\prime}=K[X_2,X_3,...,X_d]_{(X_2,X_3,...,X_d)}[X_1]$. Since the above transformation takes the maximal ideal $(\phi(X_1),X_2,...,X_d)$ of $K[X_1,X_2,...,X_d]$ to itself, then, the polynomial $h\in L^{\prime}$ is a Weierstrass polynomial. Hence  by Lemma \ref{0.6}, we have $L^{\prime}\hookrightarrow L$ is analytic isomorphism along $h$.

Hence we have $L^{\prime}[It]\hookrightarrow L[It]$ is analytic isomorphism along $h$. Since $\beta_h\in E_4(L_h[It])$, then again by applying Lemma \ref{1.37}, we have, there exists $\delta \in E_4(L[It])$  and $\theta \in SL_4(L^{\prime}[It])$ such that $\beta=\delta\theta$ and $\theta_h\in E_4(L^{\prime}_h[It])$. Since $L^{\prime}[It]= (k[X_2,X_3,...,X_d]_{(X_2,X_3,...,X_d)}[It])[X_1]$ and $h$ is monic polynomial in $X_1$ of $L^{\prime}[It]$, then by Proposition \ref{3.4}, we have, $\theta\in E_4(L^{\prime}[It])$.

Hence $\beta\in E_4(L[It])\subseteq E_4(R[It])$. Finally we have $\alpha \in E_4(R[It])$.

For the arbitrary ground field, we follow the following treatment suggested by Swan, given in \cite{lindel}.

Let $R= C_{\mathfrak{p}}$ where $C= k[X_1,X_2,\dots,X_m]/(f_1,f_2,\dots f_t)$, $\mathfrak{p}$ is a prime ideal of $C$. Let $k_0$ be the prime subfield of $k$. Choose a field extension $K$ of $k_0$ such that $K$ is finitely generated over $k_0$ and $K$ contains all coefficients of the $f_i$'s and all the elements of $k$ such that $\alpha$ is defined over $K$. Set $B= K[X_1,X_2,\dots X_m]/(f_1,f_2,\dots,f_t)$, $\mathfrak{q}= \mathfrak{p}\cap B$, $R^\prime= B_{\mathfrak{q}}$. We have $B\hookrightarrow B\otimes_Kk=C$. Now note that $R$ is a flat extension of $C$; $C$ is a flat(free) extension of $B$. Hence $R$ is flat over $B$. Hence $R$ is flat over $R^\prime$. Since $R^\prime$ is local ring, then $R$ is faithfully flat over $R^\prime$. Moreover $R^\prime\subset R$ is a local extension, i.e, $\mathfrak{m}_{R^\prime}\subset \mathfrak{m}_R$. Hence $R^{\prime}$ is a regular $K$-spot by Proposition \ref{3.10}. Hence $R^\prime$ is a regular local algebra with a separating ground field $K$ by Lemma \ref{1.39}. Hence we have $\alpha \in SL_4(R^\prime[It])\cap E_5(R^\prime[It])$ with $R^\prime$ is a regular local algebra with a separating ground field. This reduces to the case we have already considered.
$~~~~~~~~~~~~~~~~~~~~~~~~~~~~~~~~~~~~~~~~~~~\qedwhite$

\begin{theo}
Let $R$ be a regular $k$-spot of dimension $d$, $d\leq 3$, and $I\subset R$ be an ideal. Then The Vaserstein symbol $V_{R[It,t^{-1}]}:MSE_3(R[It,t^{-1}])\rightarrow W_E(R[It,t^{-1}])$ is injective.
\end{theo}

${\pf}$

By Lemma \ref{1.43}, it is enough to check that $SL_4(R[It,t^{-1}])\cap E(R[It,t^{-1}])=E_4(R[It,t^{-1}])$. Let $\alpha\in SL_4(R[It,t^{-1}])\cap E(R[It.t^{-1}])$. First assume that $R$ is a regular local algebra with a separating ground field $K$.  Let $S=R\setminus \{0\}$. By stability we have $\alpha_S\in E_4((S^{-1}R)[It,t^{-1}])$ since $\dim{(S^{-1}R)[It,t^{-1}]}\leq 1$. Hence there exists $g\in R$ such that $\alpha_h\in E_4(R_g[It,t^{-1}])$. We may assume that $g\in \mathfrak{m}$, the maximal ideal of $R$. Further we may assume that $g\in \mathfrak{m}^2$. By Theorem \ref{1.36}, there exists a subring $L$ of $R$ and an element $h\in L\cap gR$ such that $L= K[X_1,X_2,\dots,X_d]_{(\varphi(X_1),X_2,\dots,X_d)}$, where $\varphi(X_1)$ is an irreducible monic polynomial, and $L\hookrightarrow R$ is an analytic isomorphism along $h$. By Lemma \ref{2.9}, $L[It]\hookrightarrow R[It,t^{-1}]$ is an analytic isomorphism along $h$. Thus we have a patching diagram
\[
\begin{tikzcd}
L[It,t^{-1}]\arrow{r}{}\arrow{d}{} &R[It,t^{-1}]\arrow{d}{}\\
L_h[It,t^{-1}] \arrow{r}{} &R_h[It,t^{-1}]
\end{tikzcd}
\]
By Theorem \ref{1.36}, we have $g$ and $h$ are differed by a unit in $R$. Hence we have $\alpha_h\in E_4(R_h[It,t^{-1}])$. Hence by Lemma \ref{1.37} we have there exists $\beta\in SL_4(L[It,t^{-1}])$ with $\beta_h\in E_4(L_h[It,t^{-1}])$ and $\gamma \in E_4(R[It,t^{-1}])$ such that $\alpha= \gamma\beta$.  Therefore it is enough to show that $\beta \in E_4(L[It,t^{-1}])$.
 
 We may assume that $h$ belongs to the maximal ideal of $L$. Also by multiplying unit of $L$ we may assume that $h\in (\phi(X_1),X_2,...,X_d)$. Now by Lemma \ref{3.5}, we have, there is a transformation of $K[X_1,X_2,...,X_d]$, namely $X_1\mapsto X_1, X_i\mapsto X_i+\phi(X_1)^{r_i}, i\geq 2$, for some $r_i$ such that $h$ becomes a monic polynomial in $X_1$ over $K[X_2,X_3,...,X_d]$.

Let $L^{\prime}=K[X_2,X_3,...,X_d]_{(X_2,X_3,...,X_d)}[X_1]$. Since the above transformation takes the maximal ideal $(\phi(X_1),X_2,...,X_d)$ of $K[X_1,X_2,...,X_d]$ to itself, then, the polynomial $h\in L^{\prime}$ is a Weierstrass polynomial. Hence  by Lemma \ref{0.6}, we have $L^{\prime}\hookrightarrow L$ is analytic isomorphism along $h$.

Hence we have $L^{\prime}[It,t^{-1}]\hookrightarrow L[It,t^{-1}]$ is analytic isomorphism along $h$. Since $\beta_h\in E_4(L_h[It,t^{-1}])$, then again by applying Lemma \ref{1.37}, we have, there exists $\delta \in E_4(L[It,t^{-1}])$  and $\theta \in SL_4(L^{\prime}[It,t^{-1}])$ such that $\beta=\delta\theta$ and $\theta_h\in E_4(L^{\prime}_h[It,t^{-1}])$. Since $L^{\prime}[It,t^{-1}]= (k[X_2,X_3,...,X_d]_{(X_2,X_3,...,X_d)}[It,t^{-1}])[X_1]$ and $h$ is monic polynomial in $X_1$ of $L^{\prime}[It,t^{-1}]$, then by Theorem \ref{3.4}, we have, $\theta\in E_4(L^{\prime}[It,t^{-1}])$.

Hence $\beta\in E_4(L[It,t^{-1}])\subseteq E_4(R[It,t^{-1}])$. Finally we have $\alpha \in E_4(R[It,t^{-1}])$.

The case for the arbitrary ground field can be reduced to the case of the separating ground field by the same argument as given at the end of the proof of Theorem \ref{4.10}.
$~~~~~~~~~~~~~~~~~~\qedwhite$

\begin{theo}\label{extdedRees}
Let $R$ be a ring of dimension $2$ and $I\subset R$ be an ideal of $R$. Then the Vaserstein symbol $V_{R[It,t^{-1}]}: MSE_3(R[It,t^{-1}])\rightarrow W_E(R[It,t^{-1}])$ is bijective.
\end{theo}

${\pf}$ Since $\dim R[It,t^{-1}]\leq 3$, we have $V_{R[It,t^{-1}]}$ is surjective. For injectivity we must show that $SL_4(R[It,t^{-1}])\cap E(R[It,t^{-1}])= E_4(R[It,t^{-1}])$. By Corollary \ref{ex-stabi}, we have the natural map

 $\phi: GL_n(R[It,t^{-1}])/E_n(R[It,t^{-1}])\rightarrow K_1(R[It,t^{-1}])$ is an isomorphism for $n\geq 4$. Hence for $n=4$, we have $SL_4(R{It,t^{-1}})\cap E_5(R[It,t^{-1}])=E_4(R[It,t^{-1}])$.
$~~~~~~~\qedwhite$

\Addresses
\end{document}